\documentclass{amsart}

\newtheorem{theorem}{Theorem}[section]
\newtheorem{corollary}[theorem]{Corollary}
\newtheorem{lemma}[theorem]{Lemma}

\begin{document}

\title[Mail your comments to wonniepark@postech.ac.kr]{Umbilic points and Real hyperquadrics}
\author[Won K. Park]{Won K. Park\thanks{%
E-mail: wonkpark@euclid.postech.ac.kr 
\newline
} \\
...........................\\ 
Any comments, suggestions, errors to\\
wonniepark@postech.ac.kr}
\address{Department of Mathematics, Postech Pohang, Korea, 790-784}
\maketitle
\tableofcontents

\begin{abstract}
There exist polynomial identities asociated to normal form, which yield an
existence and uniqueness theorem. The space of normalized real hypersurfaces
has a natural group action. Umbilic point is defined via normal form. A
nondegenerate analytic real hypersurface is locally biholomorphic to a real
hyperquadric if and only if every point of the real hypersurface is umbilic.
\end{abstract}

\addtocounter{section}{-1}

\section{\textbf{Introduction}}

An analytic real hypersurface \thinspace $M$ is said to be in Chern-Moser
normal form if $M$ is defined by the following equation near the origin: 
\[
v=\langle z,z\rangle +\sum_{\min (s,t)\geq 2}F_{st}(z,\bar{z},u) 
\]
where 
\[
\langle z,z\rangle \equiv z^{1}\overline{z^{1}}+\cdots +z^{e}\overline{z^{e}}%
-z^{e+1}\overline{z^{e+1}}-\cdots -z^{n}\overline{z^{n}} 
\]
for a positive integer $e$ in $\frac{n}{2}\leq e\leq n,$ and 
\[
F_{st}(\mu z,\nu \bar{z},u)=\mu ^{s}\nu ^{t}F_{st}(z,\bar{z},u) 
\]
for all complex numbers $\mu ,\nu ,$ and the functions $F_{22},$ $F_{23},$ $%
F_{33}$ satisfy the condition 
\[
\Delta F_{22}=\Delta ^{2}F_{23}=\Delta ^{3}F_{33}=0. 
\]
Here the operator $\Delta $ is defined as follows: 
\begin{gather*}
\Delta \equiv D_{1}\overline{D}_{1}+\cdots +D_{e}\overline{D}_{e}-D_{e+1}%
\overline{D}_{e+1}-\cdots -D_{n}\overline{D}_{n}, \\
D_{k}=\frac{\partial }{\partial z^{k}},\quad \overline{D}_{k}=\frac{\partial 
}{\partial \overline{z^{k}}},\quad k=1,\cdots ,n.
\end{gather*}

Then we have an existence theorem of a biholomorphic normalizing mapping(cf. 
\cite{CM74}, \cite{Pa1}).

\begin{theorem}[Chern-Moser]
\label{Exi}Let $M$ be an analytic real hypersurface with nondegenerate Levi
form at the origin in $\Bbb{C}^{n+1}$ defined by the following equation: 
\[
v=F(z,\bar{z},u),\text{\quad }\left. F\right| _{0}=\left. dF\right| _{0}=0.
\]
Then there exists a biholomorphic normalizing mapping $\phi $ of $M$ such
that $\phi \left( M\right) $ is in Chern-Moser normal form.
\end{theorem}

We have a uniqueness theorem of a biholomorphic normalizing mapping(cf. \cite
{CM74}, \cite{Pa1}).

\begin{theorem}[Chern-Moser]
\label{Uni}Let $M$ be the real hypersurface in Theorem \ref{Exi} and $\phi
=(f,g)$ in $\Bbb{C}^{n}\times \Bbb{C}$ be a biholomorphic normalizing
mapping of $M$ into Chern-Moser normal form such that 
\[
\left. f\right| _{0}=\left. g\right| _{0}=0.
\]
Then $\phi $ is uniquely determined by the initial value $\sigma =(C,a,\rho
,r)$ given by 
\[
C=\left. \frac{\partial f}{\partial z}\right| _{0},\quad -Ca=\left. \frac{%
\partial f}{\partial w}\right| _{0},\quad \rho =\Re \left( \left. \frac{%
\partial g}{\partial w}\right| _{0}\right) ,\quad 2\rho r=\Re \left( \left. 
\frac{\partial ^{2}g}{\partial w^{2}}\right| _{0}\right) .
\]
\end{theorem}

Let $M$ be an analytic real hypersurface defined by the equation 
\[
v=\langle z,z\rangle +\sum_{k=3}^{\infty }F_{k}(z,\overline{z},u)
\]
and $\phi $ be a (possibly formal) biholomorphic normalizing mapping with
identity initial value in $\Bbb{C}^{n}\times \Bbb{C}$ such that 
\[
\phi =\left( z+\sum_{k=2}^{\infty }f_{k}(z,w),w+\sum_{k=3}^{\infty
}g_{k}(z,w)\right) .
\]
Suppose that the real hypersurface $\phi \left( M\right) $ is defined by the
(possibly formal) equation 
\[
v=\langle z,z\rangle +\sum_{k=3}^{\infty }F_{k}^{*}(z,\overline{z},u).
\]
Then we obtain a following family of polynomial identities(cf. \cite{CM74}): 
\begin{align}
& \left. \Re \{2\langle z,f_{m-1}\left( z,w\right) \rangle +ig_{m}\left(
z,w\right) \}\right| _{w=u+i\langle z,z\rangle }  \nonumber \\
& =F_{m}\left( z,\overline{z},u\right) -F_{m}^{*}\left( z,\overline{z}%
,u\right) +R_{m}\left( z,\overline{z},u\right) \hspace{0.5in}\text{for }%
m\geq 3,  \tag*{(0.1)}  \label{iden}
\end{align}
where $R_{m}\left( z,\overline{z},u\right) $ is a polynomial of weight $m$
consisting of the functions 
\[
f_{k-1},\quad g_{k},\quad F_{k},\quad F_{k}^{*},\quad \quad k\leq m-1.
\]

In the thesis \cite{Pa0}, we studied the polynomial identity \ref{iden} in
each weight in order to investigate the process of determining the
normalizing mapping $\phi $ and the normalized real hypersurface $\phi
\left( M\right) $. Under accepting Theorem \ref{Exi} and Theorem \ref{Uni}
to be proven independently(cf. \cite{CM74}, \cite{Pa1}), we obtain the
following existence and uniqueness theorem form the polynomial identities 
\ref{iden}:

\begin{theorem}
\label{exun}There is a natural isomorphism for each $k\geq 3$ via the
polynomial identity \ref{iden} such that 
\[
\left\{ F_{l}:l\leq k\right\} \simeq \left\{ (f_{l-1},g_{l}):l\leq k\right\}
\oplus \left\{ F_{l}^{*}:l\leq k\right\} .
\]
Hence there exist two unique mappings for each $k\geq 3$ such that 
\begin{eqnarray*}
\nu  &:&\left\{ F_{l}:l\leq k\right\} \longmapsto \left\{
(f_{l-1},g_{l}):l\leq k\right\}  \\
\kappa  &:&\left\{ F_{l}:l\leq k\right\} \longmapsto \left\{ F_{l}^{*}:l\leq
k\right\} .
\end{eqnarray*}
Further, the formal biholomorphic mapping 
\[
\phi =\left( z+\sum_{k=2}^{\infty }f_{k}(z,w),w+\sum_{k=3}^{\infty
}g_{k}(z,w)\right) 
\]
is convergent so that $\phi \left( M\right) $ is an analytic real
hypersurface defined by the equation 
\[
v=\langle z,z\rangle +\sum_{k=4}^{\infty }F_{k}^{*}\left( z,\overline{z}%
,u\right) .
\]
\end{theorem}

Kruzhilin \cite{Kr} showed that some low order terms of a normalization $%
\phi $ of $M$ are equal to the low order terms of a local automorphism of a
real hyperquadric whenever the real hypersurface $M$ is already in normal
form. We present a simple proof of Kruzhilin's Lemma. Then, as its
consequence, we show that there is a natural group action on the space of
analytic real hypersurfaces in normal form by the isotropy group $H$ of a
real hyperquadric via normalizations.

E. Cartan \cite{Ca32} and Chern-Moser \cite{CM74} have introduced umbilic
points as a local CR invariant in their geometric theory and Chern-Moser
identified the condition in normal form so that, on a nondegenerate analytic
real hypersurface $M,$ a point $p\in M$ is umbilic if there is a normal
coordinate with center at $p\in M$ such that, for $\dim M=3,$%
\begin{equation}
v=\langle z,z\rangle +\sum_{\min (s,t)\geq 2,\max (s,t)\geq 4}F_{st}(z,\bar{z%
},u)  \tag*{(0.2)}  \label{24}
\end{equation}
where 
\[
F_{24}(z,\bar{z},0)=F_{42}(z,\bar{z},0)=0 
\]
and, for $\dim M\geq 5,$%
\begin{equation}
v=\langle z,z\rangle +\sum_{\min (s,t)\geq 2}F_{st}(z,\bar{z},u) 
\tag*{(0.3)}  \label{22}
\end{equation}
where 
\[
F_{22}(z,\bar{z},0)=0. 
\]
We show that this condition of umbilic points in normal coordinate may be
taken to be the definition of umbilic points. Then we shall prove our main
theorem in this article on umbilic points and real hyperquadrics.

\begin{theorem}
\label{CK}Let $M$ be an analytic real hypersurface with nondegenerate Levi
form and $U$ be an open neighborhood of a point $p\in M$ such that $U\cap M$
consists of umbilic points. Then the open subset $U\cap M$ is locally
biholomorphic to a real hyperquadric.
\end{theorem}

We may view Theorem \ref{CK} as an analytic analogue of E. Cartan's
equivalence problem of spherical real hypersurfaces, which concerns the
local existence of biholomorphic mapping of a spherical real hypersurface $M$
to a real hyperquadric(cf. \cite{Ca32}, \cite{CM74}). In the case that the
Levi form on $M$ is definite. Then Theorem \ref{CK} assures the local
existence of biholomorphic mapping of a spherical real hypersurface $M$ to a
sphere $S^{2n+1}.$

On a nondegenerate analytic real hypersurface $M,$ a point $p\in M$ is
called spherical if there exist a neighborhood $U$ of $p$ and a
biholomorphic mapping $\phi $ on $U$ such that 
\[
\phi \left( U\cap M\right) \subset Q
\]
where $Q$ is a real hyperquadric. Then Theorem \ref{CK} is just a
characterization of a spherical point $p$ by umbilic points near $p.$

\section{Existence and uniqueness theorem}

\textbf{I}. In this article, we concern local properties of real
hypersurfaces under biholomorphic mappings so that each real hypersurface
has a distinguished point as its base point, i.e., the origin in $\Bbb{C}%
^{n+1},$ and each biholomorphic mapping leaves the origin invariant unless
specified otherwise.

Once we agree to define a real hypersurface $M$ locally by an equation of
the form 
\begin{equation}
v=F\left( z,\overline{z},u\right) .  \nonumber  \label{defining equation}
\end{equation}
Then there is one-to-one correspondence between the real hypersurface $M$
and the defining equation. Hence we may identify the real hypersurface $M$
locally with the defining equation unless serious confusion. We may abuse
notations in this regard.

We shall study some consequences of Theorem \ref{Exi} and Theorem \ref{Uni}.

\begin{lemma}[Chern-Moser]
\label{Lemma}Let $M$ be a nondegenerate analytic real hypersurface defined
by 
\[
v=F(z,\overline{z},u)\equiv \sum_{k=2}^{\infty }F_{k}(z,\overline{z},u).
\]
Let $\phi =(\sum_{k}f_{k},\sum_{k}g_{k})$ be a normalization of $M$ with
initial value $(C,a,\rho ,r)$ such that $\phi \left( M\right) $ is defined
by the equation 
\[
v=\langle z,z\rangle +\sum_{k=4}^{\infty }F_{k}^{*}(z,\overline{z},u).
\]
Then there is a family of polynomial identities, for each $k\geq 3,$%
\[
\mathcal{L}_{k}\left( f_{k-1},g_{k},F_{k}^{*}\right) =\rho F_{k}(z,\overline{%
z},u)+R_{k}(z,\overline{z},u)
\]
where $R_{k}(z,\overline{z},u)$, $k\geq 3,$ is represented by a finite
linear combination of finite multiples of the functions 
\[
f_{s-1},\quad g_{s},\quad F_{s},\quad F_{s}^{*},\quad s\leq k-1,\quad 
\mathrm{and}\text{ }\mathrm{their}\text{ }\mathrm{derivatives.}
\]
\end{lemma}

\proof
By Theorem \ref{Exi}, there is a biholomorphic normalization of $M,$ $\phi
=(f,g)$. Then we obtain the following identity 
\begin{align}
& \Im g(z,u+iF(z,\overline{z},u))  \nonumber \\
=& \langle f(z,u+iF(z,\overline{z},u)),f(z,u+iF(z,\overline{z},u))\rangle  
\nonumber \\
& +F^{*}(f(z,u+iF(z,\overline{z},u)),\overline{f(z,u+iF(z,\overline{z},u))}%
,\Re g(z,u+iF(z,\overline{z},u))).  \tag*{(1.1)}  \label{identity}
\end{align}
We expand the identity with respect to weight and collect terms of the same
weight so that

\begin{enumerate}
\item[(1)]  for weight $1$ and $2,$%
\begin{align}
\Im g_{1}(z,0)& =0  \nonumber \\
\Im g_{2}(z,u+iF_{2}(z,\overline{z},u))& =\langle
f_{1}(z,0),f_{1}(z,0)\rangle   \tag*{(1.2)}  \label{weight2}
\end{align}

\item[(2)]  for weight $m\geq 3,$%
\begin{align}
& \Im g_{m}(z,u+iF_{2}(z,\overline{z},u))-\Re g_{2}(0,1)F_{m}(z,\overline{z}%
,u)  \nonumber \\
& -2\Re \langle f_{1}(z,0),f_{m-1}(z,u+iF_{2}(z,\overline{z},u))\rangle  
\nonumber \\
& -F_{m}^{*}(f_{1}(z,0),\overline{f_{1}(z,0)},\Re g_{2}(z,u+iF_{2}(z,%
\overline{z},u)))  \nonumber \\
& =R_{m}(z,\overline{z},u)  \nonumber
\end{align}
\end{enumerate}

\noindent where the polynomial $R_{m}(z,\overline{z},u)$ is given by a
finite linear combination of finite multiples of the following functions 
\[
f_{s-1},\quad g_{s},\quad F_{s},\quad F_{s}^{*},\quad s\leq m-1,\quad 
\mathrm{and}\text{ }\mathrm{their}\text{ }\mathrm{derivatives.}
\]
With the expansion 
\[
g_{2}(z,u+iF_{2}(z,\overline{z}%
,u))=g_{2}(z,0)+g_{2}(0,1)u+ig_{2}(0,1)F_{2}(z,\overline{z},u),
\]
in the equality \ref{weight2}, we obtain 
\[
\langle f_{1}(z,0),f_{1}(z,0)\rangle =\Re g_{2}(0,1)\sum_{\alpha ,\beta
}\left( \left. \frac{\partial ^{2}F}{\partial z^{\alpha }\partial \overline{z%
}^{\beta }}\right| _{0}\right) z^{\alpha }\overline{z}^{\beta }
\]
so that 
\[
g_{2}(z,w)=\Re g_{2}(0,1)\times \{1-iF_{2}(0,0,1)\}\{w-iF_{2}(z,0,0)\}.
\]
Note that $\Re g_{2}(0,1)\neq 0$ necessarily since the mapping $\phi $ is
biholomorphic.

Then, for weight $m\geq 3$, we define a mapping $\mathcal{L}_{m}$ such that 
\begin{eqnarray*}
\mathcal{L}_{m}(f_{m-1},g_{m},F_{m}^{*}) &\equiv &\Im g_{m}(z,u+iF_{2}(z,%
\overline{z},u)) \\
&&-2\Re \langle f_{1}(z,0),f_{m-1}(z,u+iF_{2}(z,\overline{z},u))\rangle \\
&&-F_{m}^{*}(f_{1}(z,0),\overline{f_{1}(z,0)},\Re g_{2}(z,u+iF_{2}(z,%
\overline{z},u))) \\
&=&\Re g_{2}(0,1)F_{m}(z,\overline{z},u)+R_{m}(z,\overline{z},u).
\end{eqnarray*}
This completes the proof.\endproof

Let $CP_{k},$ $k\geq 3,$ be a subspace of $\Bbb{C[}z,w]\equiv \Bbb{C[}%
z^{1},\cdots ,z^{n},w]$ such that 
\[
CP_{k}=\left\{ g\in \Bbb{C[}z,w]:g(\mu z,\mu ^{2}w)=\mu ^{k}g(z,w)\right\} 
\]
and $RP_{k},$ $k\geq 3,$ be a subspace of $\Bbb{R[\Re }z,\Bbb{\Im }z,\Bbb{%
\Re }w]$ such that 
\[
RP_{k}=\left\{ F\in \Bbb{R[\Re }z,\Bbb{\Im }z,\Bbb{\Re }w]:F(\mu z,\mu 
\overline{z},\mu ^{2}u)=\mu ^{k}F(z,\overline{z},u)\right\} . 
\]
We define a subspace $NP_{k}$ of $RP_{k},$ $k\geq 3,$ such that 
\begin{eqnarray*}
NP_{k} &=&\left\{ F\in RP_{k}:F(z,\overline{z},u)=\sum_{s,t\geq 2}F_{st}(z,%
\overline{z},u)\right. \\
&&\left. \qquad \text{where }\Delta F_{22}=\Delta ^{2}F_{23}=\Delta
^{3}F_{33}=0\right\} .
\end{eqnarray*}
We shall use the following notations: 
\[
\left\{ 
\begin{array}{l}
O(m)\equiv \sum_{\mid I\mid +\mid J\mid +2k+2l\geq m}O(z^{I}\bar{z}^{J}w^{k}%
\overline{w}^{l}) \\ 
O(z^{s}\overline{z}^{t})\equiv \sum_{\mid I\mid =s,\mid J\mid =t}O(z^{I}\bar{%
z}^{J})
\end{array}
\right. 
\]
where $k,l,s,t\in \Bbb{N}$, and $I,J$ are multi-indices in $\Bbb{N}^{n}.$

\begin{lemma}
\label{cm3.2}Let $\mathcal{L}_{k}:HP_{k-1}^{n}\times HP_{k}\times
NP_{k}\longrightarrow RP_{k}$ for $k\geq 3$ be the mapping defined in Lemma 
\ref{Lemma} as follows 
\[
\mathcal{L}_{k}(f,g,F^{*})\equiv \left[ \Re \{2\langle f(z,w),Cz\rangle
+ig(z,w)\}-F^{*}(Cz,\overline{Cz},\rho \Re \chi (z,w))\right] _{w=u+iF_{2}(z,%
\overline{z},u)}
\]
where $HP_{k}^{n}=HP_{k}\times \cdots \times HP_{k}$ $(n$ times$)$ and 
\[
\left\{ 
\begin{array}{l}
(f,g,F^{*})\in HP_{k-1}^{n}\times HP_{k}\times NP_{k} \\ 
Cz=f_{1}(z,0),\quad \rho =\Re g_{2}(0,1) \\ 
\chi (z,w)=\{1-iF_{2}(0,0,1)\}\{w-iF_{2}(z,0,0)\}.
\end{array}
\right. .
\]
Then $\mathcal{L}_{k}$ is surjective for $k\geq 3$ and 
\begin{eqnarray*}
\dim _{\Bbb{R}}\ker \mathcal{L}_{3}=2n \\
\dim _{\Bbb{R}}\ker \mathcal{L}_{4}=1 \\
\dim _{\Bbb{R}}\ker \mathcal{L}_{k}=0\quad \text{for }k\geq 5.
\end{eqnarray*}
\end{lemma}

\proof
Note that 
\[
\underleftarrow{\lim }_{k}\Bbb{R}[\Re z,\Im z,u]/O(k)\simeq \Bbb{R[}[\Re
z,\Im z,u]] 
\]
where $\underleftarrow{\lim }_{k}\Bbb{R}[\Re z,\Im z,u]/O(k)$ is the inverse
limit of $\Bbb{R}[\Re z,\Im z,u]/O(k),$ $k\in \Bbb{N}.$ Hence we can take
arbitrary $F_{k}(z,\overline{z},u)$ in Lemma \ref{Lemma} so that the mapping 
$\mathcal{L}_{k},$ $k\geq 3,$ is surjective. Then by Theorem \ref{Uni} the
kernel of the mapping $\mathcal{L}_{k},$ $k\geq 3,$ is parametrized by $a,r$
in the value $(C,a,\rho ,r)$ so that 
\begin{eqnarray*}
\dim _{\Bbb{R}}\ker \mathcal{L}_{3} &=&2n \\
\dim _{\Bbb{R}}\ker \mathcal{L}_{4} &=&1 \\
\dim _{\Bbb{R}}\ker \mathcal{L}_{k} &=&0\quad \text{for }k\geq 5.
\end{eqnarray*}
This completes the proof.\endproof

The polynomial identities in Lemma \ref{Lemma} and Lemma \ref{cm3.2} yields
the following existence and uniqueness theorem.

\begin{theorem}
\label{exi-uni}Let $M$ be a nondegenerate analytic real hypersurface defined
by 
\[
v=\sum_{k=2}^{\infty }F_{k}\left( z,\overline{z},u\right) .
\]
Then there is a natural isomorphism for each $k\geq 4$ via the mapping $%
\mathcal{L}_{k},$ $k\geq 3,$ such that 
\begin{eqnarray*}
\left\{ F_{l}\in RP_{l}:l\leq k\right\} \times H \\
\simeq \left\{ (f_{l-1},g_{l})\in CP_{l-1}^{n}\times CP_{l}:l\leq k\right\}
\oplus \left\{ F_{l}^{*}\in NP_{l}:l\leq k\right\} .
\end{eqnarray*}
Hence there exist unique mappings for each $k\geq 4$ such that 
\begin{eqnarray*}
\nu :\left\{ F_{l}:l\leq k\right\} \times H\longmapsto \left\{
(f_{l-1},g_{l}):l\leq k\right\}  \\
\kappa :\left\{ F_{l}:l\leq k\right\} \times H\longmapsto \left\{
F_{l}^{*}:l\leq k\right\} .
\end{eqnarray*}
Further, the formal biholomorphic mapping 
\[
\phi =\left( \sum_{k=1}^{\infty }f_{k},\sum_{k=2}^{\infty }g_{k}\right) 
\]
converges so that $\phi \left( M\right) $ is an analytic real hypersurface
defined by the equation 
\[
v=\langle z,z\rangle +\sum_{k=4}^{\infty }F_{k}^{*}\left( z,\overline{z}%
,u\right) .
\]
\end{theorem}

Hence, for the special case of $k=\infty $ in Theorem \ref{exi-uni}, we
obtain

\begin{theorem}
\label{main}Let $M$ be a nondegenerate analytic real hypersurface defined by 
\[
v=\sum_{k=2}^{\infty }F_{2}(z,\overline{z},u)
\]
and $\phi =(\sum_{k}f_{k},\sum_{k}g_{k})$ be a normalization of $M$ such
that the real hypersurface $\phi \left( M\right) $ is defined in normal form
by the equation 
\[
v=\langle z,z\rangle +\sum_{k}^{\infty }F_{k}^{*}(z,\overline{z},u).
\]
Then the functions $f_{k-1},g_{k},F_{k}^{*},$ $k\geq 3,$ are given as a
finite linear combination of finite multiples of the following factors:

\begin{enumerate}
\item[(1)]  the coefficients of the functions $F_{l}$, $l\leq k,$

\item[(2)]  the constants $C,C^{-1},\rho ,\rho ^{-1},a,r,$
\end{enumerate}

where $(C,a,\rho ,r)$ are the initial value of the normalization $\phi .$
\end{theorem}

\textbf{II}. We shall examine concrete versions of Theorem \ref{exi-uni} as
existence and uniqueness theorem. Let $h=(f,g)$ and $\phi =(\widetilde{f},%
\widetilde{g})$ be holomorphic mappings in $\Bbb{C}^{n}\times \Bbb{C}$ such
that 
\[
f=\widetilde{f}+O(k)\text{ and\textit{\ }}g=\widetilde{g}+O(k+1). 
\]
Then we shall write 
\[
h=\phi +O_{\times }(k+1). 
\]

\begin{theorem}
\label{Cor.3.2}Let $M$ be an analytic real hypersurface with nondegenerate
Levi form defined by 
\begin{equation}
v=F(z,\bar{z},u),\quad \left. F\right| _{0}=\left. dF\right| _{0}=0. 
\tag*{(1.3)}  \label{2}
\end{equation}
If $h=(f,g)$ is a biholomorphic mapping such that 
\begin{eqnarray*}
f(z,w) &=&C(z-aw)+f^{*}(z,w), \\
g(z,w) &=&\rho (w+rw^{2})+g^{*}(z,w),
\end{eqnarray*}
where the functions $f^{*}(z,w)$ and $g^{*}(z,w)$ satisfy the condition: 
\begin{equation}
\left. f^{*}\right| _{0}=\left. df^{*}\right| _{0}=\left. g^{*}\right|
_{0}=\left. dg^{*}\right| _{0}=\Re \left( \left. g_{ww}^{*}\right|
_{0}\right) =0,  \tag*{(1.4)}  \label{ini-cond}
\end{equation}
and if the transformed real hypersurface $h\left( M\right) $ is defined by 
\begin{equation}
v=F^{*}(z,\bar{z},u)+O(k+1),  \tag*{(1.5)}  \label{3.7}
\end{equation}
where $v=F^{*}(z,\bar{z},u)$ is in normal form, then there is a
normalization of $M,$ $\phi _{\sigma }$, with initial value $\sigma
=(C,a,\rho ,r)\in H$ such that 
\[
h=\phi _{\sigma }+O_{\times }(k+1).
\]
\end{theorem}

\proof
Note that 
\[
F_{m}^{*}\in NP_{m}\quad \text{for }3\leq m\leq k. 
\]
Then we obtain, for $3\leq m\leq k,$%
\[
\mathcal{L}_{m}(f_{m-1},g_{m},F_{m}^{*})=\rho F_{m}(z,\overline{z}%
,u)+R_{m}(z,\overline{z},u) 
\]
where $R_{m}(z,\overline{z},u)$ is a linear combination of multiples of $%
f_{s-1},g_{s},F_{s},F_{s}^{*},$ $s\leq m-1,$ and their derivatives. The
condition \ref{ini-cond} makes the mappings $\mathcal{L}_{k},$ $k\geq 3,$ be
one-to-one. Hence the coefficients of the functions 
\[
f_{m-1},\quad g_{m},\quad F_{m}^{*}\quad \text{for }m\leq k 
\]
are completely determined by the coefficients of the functions $F_{m},$ $%
m\leq k,$ and the initial value $\sigma =(C,a,\rho ,r).$

By Lemma \ref{Lemma} and Lemma \ref{cm3.2}, the normalization $\phi _{\sigma
}$ and the defining function $\phi _{\sigma }(M)$ is uniquely determined by
the value $\sigma =(C,a,\rho ,r)$ via the following equalities: 
\[
\mathcal{L}_{m}(f_{m-1}^{*},g_{m}^{*},F_{m}^{**})=\rho F_{m}(z,\overline{z}%
,u)+R_{m}(z,\overline{z},u) 
\]
where 
\[
\phi _{\sigma }(M):v=\langle z,z\rangle +F^{**}\left( z,\overline{z}%
,u\right) 
\]
and 
\begin{eqnarray*}
\phi _{\sigma } &=&(f^{*},g^{*}) \\
f^{*}(z,w) &=&\sum_{m=1}^{\infty }f_{m}^{*}(z,w) \\
g^{*}(z,w) &=&\sum_{m=1}^{\infty }g_{m}^{*}(z,w) \\
F^{**}(z,\overline{z},u) &=&\sum_{m=4}^{\infty }F_{m}^{**}(z,\overline{z},u).
\end{eqnarray*}
Since $\phi _{\sigma }$ and $h$ have the same initial value, we obtain 
\[
f_{m-1}=f_{m-1}^{*},\quad g_{m}=g_{m}^{*}\quad \text{for }m\leq k 
\]
so that 
\[
h=\phi _{\sigma }+O_{\times }(k+1). 
\]
This completes the proof.\endproof

\begin{lemma}
\label{basic}Let $\varphi =(f,g),$ $\varphi _{1}=(f_{1},g_{1}),$ $\varphi
_{2}=(f_{2},g_{2})$ be biholomorphic mappings in $\Bbb{C}^{n}\times \Bbb{C}$
such that 
\[
\left\{ 
\begin{array}{l}
\left. f\right| _{0}=\left. g\right| _{0}=\left. \frac{\partial g}{\partial z%
}\right| _{0}=0 \\ 
\left. f_{1}\right| _{0}=\left. g_{1}\right| _{0}=\left. \frac{\partial g_{1}%
}{\partial z}\right| _{0}=0 \\ 
\left. f_{2}\right| _{0}=\left. g_{2}\right| _{0}=\left. \frac{\partial g_{2}%
}{\partial z}\right| _{0}=0.
\end{array}
\right. 
\]
Then for each $k\geq 3$

\begin{enumerate}
\item[(1)]  $\varphi ^{-1}=\phi ^{-1}+O_{\times }(k)$ whenever 
\[
\varphi =\phi +O_{\times }(k),
\]

\item[(2)]  $\varphi _{1}\circ \varphi _{2}=\phi _{1}\circ \phi
_{2}+O_{\times }(k)$ whenever 
\begin{eqnarray*}
\varphi _{1}=\phi _{1}+O_{\times }(k) \\
\varphi _{2}=\phi _{2}+O_{\times }(k).
\end{eqnarray*}
\end{enumerate}
\end{lemma}

\proof
Let $\varphi =(f,g)$ and $\varphi ^{-1}=(f^{*},g^{*})$ so that we obtain
following identities 
\begin{equation}
\left\{ 
\begin{array}{l}
f\left( f^{*}\left( z,w\right) ,g^{*}\left( z,w\right) \right) =z \\ 
g\left( f^{*}\left( z,w\right) ,g^{*}\left( z,w\right) \right) =w
\end{array}
\right. .  \tag*{(1.6)}  \label{inverse}
\end{equation}
We expand the identity \ref{inverse} and collect terms of the same weight
with the weight decompositions 
\[
\left\{ 
\begin{array}{l}
f\left( z,w\right) =\sum_{k=1}^{\infty }f_{k}\left( z,w\right) ,\quad
g\left( z,w\right) =\sum_{k=2}^{\infty }g_{k}\left( z,w\right)  \\ 
f^{*}\left( z,w\right) =\sum_{k=1}^{\infty }f_{k}^{*}\left( z,w\right)
,\quad g^{*}\left( z,w\right) =\sum_{k=2}^{\infty }g_{k}^{*}\left(
z,w\right) 
\end{array}
\right. .
\]
We set 
\[
C=\left( \left. \frac{\partial f}{\partial z}\right| _{0}\right) \quad 
\mathrm{and}\quad \rho =\left( \left. \frac{\partial g}{\partial w}\right|
_{0}\right) .
\]
Then we easily see that the function 
\[
Cf_{m-1}^{*}\left( z,w\right) 
\]
is given by a finite linear combination of finite multiples of the functions 
\[
f_{1},\cdots ,f_{m-1},f_{1}^{*},\cdots ,f_{m-2}^{*},g_{2}^{*},\cdots
,g_{m-1}^{*}\quad \text{and their derivatives,}
\]
and the function 
\[
\rho g_{m}^{*}\left( z,w\right) 
\]
is given by a finite linear combination of finite multiples of the functions 
\[
g_{2},\cdots ,g_{m},f_{1}^{*},\cdots ,f_{\left[ \frac{m}{2}\right]
}^{*},g_{2}^{*},\cdots ,g_{m-1}^{*}\quad \text{and their derivatives.}
\]
Then, by a simple induction argument, the functions 
\[
\left( f_{m-1}^{*}\left( z,w\right) ,g_{m}^{*}\left( z,w\right) \right) 
\]
are given by a finite linear combination of finite multiples of the
functions 
\[
f_{1},\cdots ,f_{m-1},g_{2},\cdots ,g_{m}\quad \text{and their derivatives}
\]
and 
\[
C,\quad C^{-1},\quad \rho ,\quad \rho ^{-1}.
\]
Hence we have proved 
\begin{equation}
\varphi ^{-1}=\phi ^{-1}+O_{\times }(k)  \tag*{(1.7)}  \label{resu}
\end{equation}
whenever 
\[
\varphi =\phi +O_{\times }(k).
\]

Since the mapping $\varphi _{1}=(f_{1},g_{1})$ satisfies the condition 
\[
\left. f_{1}\right| _{0}=\left. g_{1}\right| _{0}=\left. \frac{\partial g_{1}%
}{\partial z}\right| _{0}=0, 
\]
we obtain 
\[
\varphi _{1}\circ \varphi _{2}=\varphi _{1}\circ \phi _{2}+O_{\times }(k) 
\]
whenever 
\[
\varphi _{2}=\phi _{2}+O_{\times }(k). 
\]
Then, by the result \ref{resu}, we have 
\begin{equation}
\varphi _{2}^{-1}\circ \varphi _{1}^{-1}=\phi _{2}^{-1}\circ \varphi
_{1}^{-1}+O_{\times }(k).  \tag*{(1.8)}  \label{step1}
\end{equation}
Note that the mapping $\phi _{2}^{-1}=(f_{2}^{*},g_{2}^{*})$ satisfies the
condition 
\begin{equation}
\left. f_{2}^{*}\right| _{0}=\left. g_{2}^{*}\right| _{0}=\left. \frac{%
\partial g_{2}^{*}}{\partial z}\right| _{0}=0,  \tag*{(1.9)}  \label{condi}
\end{equation}
whenever $\phi _{2}=(f_{2},g_{2})$ satisfies 
\[
\left. f_{2}\right| _{0}=\left. g_{2}\right| _{0}=\left. \frac{\partial g_{2}%
}{\partial z}\right| _{0}=0. 
\]
By the condition \ref{condi}, we obtain 
\begin{equation}
\phi _{2}^{-1}\circ \varphi _{1}^{-1}=\phi _{2}^{-1}\circ \phi
_{1}^{-1}+O_{\times }(k)  \tag*{(1.10)}  \label{step2}
\end{equation}
whenever 
\[
\varphi _{1}^{-1}=\phi _{1}^{-1}+O_{\times }(k). 
\]
By the result \ref{resu}, the equalities \ref{step1} and \ref{step2} yields 
\begin{eqnarray*}
\varphi _{2}^{-1}\circ \varphi _{1}^{-1} &=&\phi _{2}^{-1}\circ \varphi
_{1}^{-1}+O_{\times }(k) \\
&=&\phi _{2}^{-1}\circ \phi _{1}^{-1}+O_{\times }(k).
\end{eqnarray*}
Once more by the result \ref{resu}, we obtain 
\[
\varphi _{1}\circ \varphi _{2}=\phi _{1}\circ \phi _{2}+O_{\times }(k). 
\]
This completes the proof.\endproof

\begin{lemma}
Let $M$ be a real hypersurface defined by 
\begin{eqnarray*}
v=F(z,\bar{z},u), \\
\left. F\right| _{0}=\left. F_{z}\right| _{0}=\left. F_{\overline{z}}\right|
_{0}=0.
\end{eqnarray*}
Let $\phi _{1}=(f,g),\phi _{2}=(f^{\prime },g^{\prime })$ be biholomorphic
mappings of $M$ such that 
\begin{eqnarray*}
\left. f\right| _{0}=\left. g\right| _{0}=\left. g_{z}\right| _{0}=0 \\
\left. f^{\prime }\right| _{0}=\left. g^{\prime }\right| _{0}=\left.
g_{z}^{\prime }\right| _{0}=0
\end{eqnarray*}
and, for $k\geq 3,$%
\[
\phi _{1}=\phi _{2}+O_{\times }(k+1).
\]
Suppose that the transformed real hypersurfaces $\phi _{1}\left( M\right) ,$ 
$\phi _{2}\left( M\right) $ are defined by 
\[
v=G\left( z,\overline{z},u\right) ,\quad v=G^{\prime }\left( z,\overline{z}%
,u\right) .
\]
Then 
\[
G\left( z,\overline{z},u\right) =G^{\prime }\left( z,\overline{z},u\right)
+O(k+1).
\]
\end{lemma}

\proof
Note that 
\[
g(z,w)=g^{\prime }(z,w)+O(k+1) 
\]
and 
\[
F(z,\overline{z},u)=O(2) 
\]
which yields 
\[
g(z,u+iF(z,\overline{z},u))=g^{\prime }(z,u+iF(z,\overline{z},u))+O(k+1). 
\]
Then we obtain 
\begin{eqnarray*}
G(z^{*},\overline{z}^{*},u^{*}) &=&\frac{1}{2i}\left. \left\{ g(z,w)-%
\overline{g}(\overline{z},\overline{w}))\right\} \right| _{(z,w)=\phi
_{1}^{-1}(z^{*},w^{*})} \\
&=&\frac{1}{2i}\left. \left\{ g^{\prime }(z,w)-\overline{g^{\prime }}(%
\overline{z},\overline{w}))\right\} \right| _{(z,w)=\phi
_{1}^{-1}(z^{*},w^{*})}+O(k+1)
\end{eqnarray*}
By Lemma \ref{basic}, we have 
\[
\phi _{1}^{-1}=\phi _{2}^{-1}+O_{\times }(k+1) 
\]
and note that 
\[
\left. g_{z}\right| _{0}=\left. g_{z}^{\prime }\right| _{0}=0. 
\]
Then we obtain 
\begin{eqnarray*}
G(z^{*},\overline{z}^{*},u^{*}) &=&\frac{1}{2i}\left. \left\{ g^{\prime
}(z,w)-\overline{g^{\prime }}(\overline{z},\overline{w}))\right\} \right|
_{(z,w)=\phi _{1}^{-1}(z^{*},w^{*})}+O(k+1) \\
&=&\frac{1}{2i}\left. \left\{ g^{\prime }(z,w)-\overline{g^{\prime }}(%
\overline{z},\overline{w}))\right\} \right| _{(z,w)=\phi
_{2}^{-1}(z^{*},w^{*})}+O(k+1).
\end{eqnarray*}
Note that 
\[
G^{\prime }(z^{*},\overline{z}^{*},u^{*})=\frac{1}{2i}\left. \left\{
g^{\prime }(z,w)-\overline{g^{\prime }}(\overline{z},\overline{w}))\right\}
\right| _{(z,w)=\phi _{2}^{-1}(z^{*},w^{*})}. 
\]
Thus we obtain 
\[
G(z^{*},\overline{z}^{*},u^{*})=G^{\prime }(z^{*},\overline{z}%
^{*},u^{*})+O(k+1). 
\]
This completes the proof.\endproof

Hence we obtain the following theorem

\begin{theorem}
\label{Th3}Let $M$ be an analytic real hypersurface in Theorem \ref{Cor.3.2}
and $h=(f,g)$ be a biholomorphic mapping such that 
\[
h=\phi +O_{\times }(k+1)
\]
where $\phi $ is a normalization of $M$ with initial value $\sigma
=(C,a,\rho ,r)$. Suppose that 
\begin{eqnarray*}
h\left( M\right) :v=G\left( z,\overline{z},u\right)  \\
\phi \left( M\right) :v=G^{\prime }\left( z,\overline{z},u\right) 
\end{eqnarray*}
Then 
\[
G\left( z,\overline{z},u\right) =G^{\prime }\left( z,\overline{z},u\right)
+O(k+1).
\]
\end{theorem}

\textbf{III}. Let $\phi $ be a fractional linear mapping such that 
\begin{equation}
\phi =\phi _{\sigma }:\left\{ 
\begin{array}{c}
z^{*}=\frac{C(z-aw)}{1+2i\langle z,a\rangle -w(r+i\langle a,a\rangle )} \\ 
w^{*}=\frac{\rho w}{1+2i\langle z,a\rangle -w(r+i\langle a,a\rangle )}
\end{array}
\right.  \tag*{(1.11)}  \label{2.1}
\end{equation}
where the constants $\sigma =(C,a,\rho ,r)$ satisfy 
\begin{gather*}
a\in \Bbb{C}^{n},\quad \rho \neq 0,\quad \rho ,r\in \Bbb{R}, \\
C\in GL(n;\Bbb{C}),\quad \langle Cz,Cz\rangle =\rho \langle z,z\rangle .
\end{gather*}
Note that the mapping $\phi $ decomposes to 
\[
\phi =\varphi \circ \psi , 
\]
where 
\begin{equation}
\psi :\left\{ 
\begin{array}{c}
z^{*}=\frac{z-aw}{1+2i\langle z,a\rangle -i\langle a,a\rangle w} \\ 
w^{*}=\frac{w}{1+2i\langle z,a\rangle -i\langle a,a\rangle w}
\end{array}
\right. \quad \text{and}\quad \varphi :\left\{ 
\begin{array}{c}
z^{*}=\frac{Cz}{1-rw} \\ 
w^{*}=\frac{\rho w}{1-rw}
\end{array}
\right. .  \tag*{(1.12)}  \label{2.2}
\end{equation}
We easily verify 
\[
\phi ^{*}(v-\langle z,z\rangle )=(v-\langle z,z\rangle )\rho (1+\delta
)^{-1}(1+\overline{\delta })^{-1}, 
\]
where 
\[
1+\delta =1+2i\langle z,a\rangle -(r+i\langle a,a\rangle )w. 
\]
By Theorem \ref{Uni}, each element of the isotropy subgroup of the
automorphism group of a real hyperquadric $v=\langle z,z\rangle $ is given
by a fractional linear mapping $\phi _{\sigma }$ in \ref{2.1}.

\begin{theorem}
\label{Cor.4.6}Let $M$ be an analytic real hypersurface in normal form such
that 
\[
v=\langle z,z\rangle +F_{l}(z,\bar{z},u)+\sum_{k\geq l+1}F_{k}(z,\bar{z},u),
\]
where 
\[
F_{l}(z,\bar{z},u)\neq 0.
\]
Let $N_{\sigma }$ be a normalization of $M$ and $\phi _{\sigma }$ be an
automorphism of the real hyperquadric with the initial value $\sigma
=(C,a,\rho ,r)\in H$. Suppose that the transformed real hypersurface $%
N_{\sigma }(M)$ is defined by 
\[
v=\langle z,z\rangle +F^{*}(z,\bar{z},u).
\]
Then

\begin{enumerate}
\item[(1)]  $N_{\sigma }=\phi _{\sigma }+O_{\times }(l+1),$

\item[(2)]  $F^{*}(z,\bar{z},u)=\rho F_{l}(C^{-1}z,\overline{C^{-1}z},\rho
^{-1}u)+O(l+1).$
\end{enumerate}
\end{theorem}

\proof
We easily compute $\psi ^{-1}$ and $\varphi ^{-1}$ as follows: 
\[
\psi ^{-1}:\left\{ 
\begin{array}{c}
z=\frac{z^{*}+aw^{*}}{1-2i\langle z^{*},a\rangle -i\langle a,a\rangle w^{*}}
\\ 
w=\frac{w^{*}}{1-2i\langle z^{*},a\rangle -i\langle a,a\rangle w^{*}}
\end{array}
\right. \quad \text{\textit{and}}\quad \varphi ^{-1}:\left\{ 
\begin{array}{c}
z=\frac{C^{-1}z^{*}}{1+r\rho ^{-1}w^{*}} \\ 
w=\frac{\rho ^{-1}w^{*}}{1+r\rho ^{-1}w^{*}}
\end{array}
\right. . 
\]
Since $\phi _{\sigma }^{-1}=\psi ^{-1}\circ \varphi ^{-1},$ we obtain 
\[
\phi _{\sigma }^{-1}=\phi _{\sigma ^{-1}}:\left\{ 
\begin{array}{c}
z=\frac{C^{-1}(z^{*}+\rho ^{-1}Caw^{*})}{1-2i\langle z^{*},\rho
^{-1}Ca\rangle -w^{*}(-r\rho ^{-1}+i\langle \rho ^{-1}Ca,\rho ^{-1}Ca\rangle
)} \\ 
w=\frac{\rho ^{-1}w^{*}}{1-2i\langle z^{*},\rho ^{-1}Ca\rangle -w^{*}(-r\rho
^{-1}+i\langle \rho ^{-1}Ca,\rho ^{-1}Ca\rangle )}
\end{array}
\right. 
\]
where 
\[
\sigma ^{-1}=(C^{-1},-\rho ^{-1}Ca,\rho ^{-1},-r\rho ^{-1})\in H. 
\]
Thus we have 
\[
v-\langle z,z\rangle =(v^{*}-\langle z^{*},z^{*}\rangle )\rho ^{-1}(1-\delta
^{*})^{-1}(1-\overline{\delta ^{*}})^{-1}, 
\]
where 
\[
1-\delta ^{*}=1-2i\rho ^{-1}\langle z^{*},Ca\rangle -\rho
^{-1}w^{*}(-r+i\langle a,a\rangle ). 
\]

By the mapping $\phi _{\sigma }$ in the decomposition $N_{\sigma }=E\circ
\phi _{\sigma },$ we obtain that 
\begin{eqnarray*}
v^{*} &=&\langle z^{*},z^{*}\rangle +\rho F_{l}(C^{-1}z^{*},\overline{%
C^{-1}z^{*}},\rho ^{-1}u^{*})+\sum_{\mid I\mid +\mid J\mid +2k\geq
l+1}O(z^{*I}\bar{z}^{*J}u^{*k}), \\
&=&\langle z^{*},z^{*}\rangle +\rho F_{l}(C^{-1}z^{*},\overline{C^{-1}z^{*}}%
,\rho ^{-1}u^{*})+O(l+1).
\end{eqnarray*}
Note that the following real hypersurface is in normal form: 
\[
v=\langle z,z\rangle +\rho F_{l}(C^{-1}z,\overline{C^{-1}z},\rho ^{-1}u). 
\]
By Theorem \ref{Cor.3.2}, $N_{\sigma }$ agrees with $\phi _{\sigma }$ up to
weight $l$ so that 
\[
N_{\sigma }=\phi _{\sigma }+O_{\times }(l+1). 
\]
Then Theorem \ref{Th3} implies 
\[
F^{*}\left( z,\overline{z},u\right) =\rho F_{l}(C^{-1}z,\overline{C^{-1}z}%
,\rho ^{-1}u)+O(l+1). 
\]
This completes the proof.\endproof

Since $\det (C)\neq 0$ and $\rho \neq 0$, the weight $l$ of the real
hypersurface $M$ in Theorem \ref{Cor.4.6} is the non-vanishing lowest weight
of the function $F^{*}\left( z,\overline{z},u\right) $ regardless of the
initial value $\sigma $, where 
\[
N_{\sigma }(M):v=\langle z,z\rangle +F^{*}\left( z,\overline{z},u\right) . 
\]
Thus the non-vanishing lowest weight $l$ is an invariant of $M$ under germs
of biholomorphic mappings.

\begin{corollary}
\label{Co1}Let $M$ be an analytic real hypersurface in normal form which is
invariant under the action of all normalizations. Then necessarily, $M$ is a
real hyperquadric.
\end{corollary}

\proof
Suppose that $M$ is not a real hyperquadric defined by 
\[
v=\langle z,z\rangle +F\left( z,\overline{z},u\right) . 
\]
Then there is a positive integer $l$ such that 
\[
F\left( z,\overline{z},u\right) =F_{l}(z,\overline{z},u)+O(l+1), 
\]
where 
\[
F_{l}(z,\overline{z},u)\neq 0. 
\]
We take an initial value $\sigma =(C,a,\rho ,r)$ such that 
\[
C=\sqrt{\rho },\text{ }a=0,\text{ }\rho >0,\text{ }r=0. 
\]
Let $\phi $ be a normalization of $M$ with initial value $\sigma =(C,a,\rho
,r)$ such that 
\[
\phi \left( M\right) :v=\langle z,z\rangle +F^{*}\left( z,\overline{z}%
,u\right) 
\]
Then by Theorem \ref{Cor.4.6}, we obtain 
\begin{eqnarray*}
F^{*}\left( z,\overline{z},u\right) &=&\rho F_{l}(C^{-1}z,\overline{C^{-1}z}%
,\rho ^{-1}u)+O(l+1) \\
&=&\rho ^{\frac{2-l}{2}}F_{l}(z,\overline{z},u)+O(l+1),
\end{eqnarray*}
where 
\[
\rho >0,\quad \sigma =(\sqrt{\rho },0,\rho ,0)\in H. 
\]
Since $l\geq 4$ necessarily$,$ the assumption implies, for all $\sigma \in
H, $%
\[
F^{*}\left( z,\overline{z},u\right) =F_{l}(z,\overline{z},u)+O(l+1). 
\]
Hence we obtain 
\[
F_{l}(z,\overline{z},u)=0. 
\]
This is a contradiction to the choice of the integer $l.$ Thus we have
showed $F\left( z,\overline{z},u\right) =0$ so that $M$ is a real
hyperquadric.\endproof

\begin{lemma}[Kruzhilin]
\label{Lem.4.1}Let $M$ be an analytic real hypersurface in normal form. Then
each normalization $N_{\sigma }$ with initial value $\sigma \in H$ satisfies
the following relation: 
\[
N_{\sigma }=\phi _{\sigma }+O_{\times }(5)
\]
where $\phi _{\sigma }$ is a local automorphism of a real hyperquadric with
initial value $\sigma \in H$.
\end{lemma}

\proof
Note that $M$ is of the form: 
\[
v=\langle z,z\rangle +F_{22}(z,\bar{z},0)+O(5), 
\]
where $F_{22}(z,\bar{z},0)$ is of weight 4. By Theorem \ref{Cor.4.6}, the
normalization $N_{\sigma }$ agrees with the mapping $\phi _{\sigma }$ in low
order terms such that 
\[
N_{\sigma }=\phi _{\sigma }+O_{\times }(5). 
\]
This completes the proof.\endproof

\begin{theorem}
\label{Cor.4.7}Let $M$ be an analytic real hypersurface in normal form and $%
N_{\sigma _{2}}$ be a normalization of $M$ with initial value $\sigma
_{2}\in H$. Let $N_{\sigma _{1}}$ be a normalization of $N_{\sigma
_{2}}\left( M\right) $ with initial value $\sigma _{1}\in H$ and $N_{\sigma
_{1}\sigma _{2}}$ be a normalization of $M$ with initial value $\sigma
_{1}\sigma _{2}\in H$. Then 
\[
N_{\sigma _{1}}\circ N_{\sigma _{2}}=N_{\sigma _{1}\sigma _{2}}.
\]
\end{theorem}

\proof
Since $M$ is in normal form, Lemma \ref{Lem.4.1} yields 
\[
N_{\sigma }=\phi _{\sigma }+O_{\times }(5) 
\]
for all $\sigma \in H.$ Hence we obtain 
\begin{eqnarray*}
N_{\sigma _{1}}\circ N_{\sigma _{2}} &=&\phi _{\sigma _{1}}\circ \phi
_{\sigma _{2}}+O_{\times }(5) \\
&=&\phi _{\sigma _{1}\sigma _{2}}+O_{\times }(5) \\
&=&N_{\sigma _{1}\sigma _{2}}+O_{\times }(5).
\end{eqnarray*}
Note that the initial value $\sigma $ of a normalization $N_{\sigma }=(f,g)$
is completely determined by the terms 
\begin{eqnarray*}
&&f_{1}(z,w),\quad f_{2}(z,w),\quad f_{3}(z,w) \\
&&g_{2}(z,w),\quad g_{3}(z,w),\quad g_{4}(z,w)
\end{eqnarray*}
where 
\[
f(z,w)=\sum_{k=1}^{\infty }f_{k}(z,w),\quad g(z,w)=\sum_{k=2}^{\infty
}g_{k}(z,w). 
\]
Then Theorem \ref{Uni} yields 
\[
N_{\sigma _{1}}\circ N_{\sigma _{2}}=N_{\sigma _{1}\sigma _{2}}. 
\]
This completes the proof.\endproof

\textbf{IV}. For a real hypersurface $M$ in normal form, we define the
isotropy subgroup $H(M)$ of $M$ as follows: 
\[
H(M)=\{\sigma \in H:N_{\sigma }(M)=M\}. 
\]
It is known that $H(M)$ is a Lie group(cf. \cite{Pa1}).

\begin{lemma}
\label{Thm.5.7}Let $M$ be a nondegenerate analytic real hypersurface defined
by 
\[
v=F(z,\overline{z},u),\quad \left. F\right| _{0}=\left. dF\right| _{0}=0
\]
and $M^{\prime }$ be another analytic real hypersurface in normal form which
is biholomorphic to $M$ near the origin. Then there is an element $\sigma
\in H$ and a normalization $N_{\sigma }$ such that $M^{\prime }=N_{\sigma
}(M).$ Suppose in addition that $M$ is in normal form. Then 
\[
M^{\prime }=N_{\sigma ^{\prime }}(M)\quad \text{if and only if}\quad \sigma
^{\prime }\in H(M^{\prime })\sigma =\sigma H(M),
\]
where $H(M^{\prime })\sigma $ is a left coset of $H(M^{\prime })$ in $H$ and 
$\sigma H(M)$ is a right coset of $H(M)$ in $H.$
\end{lemma}

\proof
Since $M,M^{\prime }$ are biholomorphic, we take a biholomorphic mapping $%
N=(f,g).$ Then by Theorem \ref{Uni} we have $N=N_{\sigma }$ with $\sigma
=(C,a,\rho ,r),$ where 
\[
C=\left. \frac{\partial f}{\partial z}\right| _{0},\quad -Ca=\left. \frac{%
\partial f}{\partial w}\right| _{0},\quad \rho =\left( \left. \frac{\partial
g}{\partial w}\right| _{0}\right) ,\quad 2\rho r=\Re \left( \left. \frac{%
\partial ^{2}g}{\partial w^{2}}\right| _{0}\right) . 
\]
Suppose that $M$ is in normal form. Then by Theorem \ref{Cor.4.7} we have 
\begin{eqnarray*}
M^{\prime } &=&N_{\sigma ^{\prime }}(M)=N_{\sigma ^{\prime }\sigma
^{-1}}\circ N_{\sigma }(M) \\
&=&N_{\sigma ^{\prime }\sigma ^{-1}}(M^{\prime }),
\end{eqnarray*}
which yields $\sigma ^{\prime }\sigma ^{-1}\in H(M^{\prime }).$ Hence we
obtain 
\[
\sigma ^{\prime }\in H(M^{\prime })\sigma . 
\]
Clearly, the converse is also true.

Suppose that $M,M^{\prime }$ are both in normal form. Then $M^{\prime
}=N_{\sigma }(M)$ if and only if $M=N_{\sigma ^{-1}}(M^{\prime }).$ Thus we
obtain 
\[
\sigma ^{-1}H(M^{\prime })\sigma \subset H(M),\quad \sigma H(M)\sigma
^{-1}\subset H(M^{\prime }) 
\]
so that 
\[
\sigma H(M)=H(M^{\prime })\sigma 
\]
This completes the proof.\endproof

\begin{theorem}
Let $\frak{N}$ be the set of real hypersurfaces in normal form with
signature $(s,n-s),$ $\frac{n}{2}\leq s\leq n.$ Then $\frak{N}$ is an $H$%
-space under the group $H$ action via normalizations.
\end{theorem}

\proof
Let $M$ be a real hypersurface in normal form. Then, by Theorem \ref{Uni},
the normalization of $M$ with identity initial value is necessarily the
identity map. Further, by Theorem \ref{Cor.4.7}, we have 
\[
N_{\sigma _{1}}\circ N_{\sigma _{2}}(M)=N_{\sigma _{1}\sigma _{2}}(M). 
\]
Thus the set $\frak{N}$ is an $H$-space under the group $H$ action via
normalizations. This completes the proof.\endproof

Let $\frak{N}/H$ be the orbit space of the $H$-space $\frak{N}.$ Suppose
that $M,M^{\prime \prime }$ be analytic real hypersurfaces in normal form
which are biholomorphic (even formally) to each other near the origin. Then,
by Lemma \ref{Thm.5.7}, there is a normalization $N_{\sigma }$ with an
initial value $\sigma \in H$ such that 
\[
M^{\prime \prime }=N_{\sigma }(M). 
\]
Hence the orbit space $\frak{N}/H$ may be isomorphic to the orbit space $%
\frak{M}$ of germs of strictly pseudoconvex analytic real hypersurfaces
under germs of biholomorphic mappings.

We have already seen several local invariants under germs of biholomorphic
mappings such as the signature $(s,n-s)$ of Levi form on $M$ and the
nonvanishing lowest weight $l$ of the defining function of $M$ in a normal
coordinate$.$ It is known that every CR diffeomorphism of class $C^{1}$
between analytic real hypersurfaces with nondegenerate Levi form is
necessarily real-analytic so as to extend to a biholomorphic mapping on a
neighborhood(cf. \cite{Pi75}, \cite{Le77}, \cite{BJT85}). Therefore, we may
conclude that the orbit space $\frak{N}/H$ is a classifying space of germs
of nondegenerate analytic real hypersurfaces. Hence any function defined via
Chern-Moser normal coordinate is a local CR invariant whenever it does not
depend on the choice of normal coordinates. We refer to Wells \cite{Wells}\
for a history of many aspects of Cauchy-Riemann invariants.

\section{Umbilic points on analytic real hypersurfaces}

\textbf{I}. Let $M$ be an analytic real hypersurface and $p$ be a point on $%
M $ such that $M$ is, in normal coordinate with center at $p,$ defined by 
\[
v=\langle z,z\rangle +\sum_{\min (s,t)\geq 2}F_{st}(z,\bar{z},u), 
\]
where 
\[
\Delta F_{22}=\Delta ^{2}F_{23}=\Delta ^{3}F_{33}=0. 
\]
Note that $F_{22}=F_{23}=F_{33}=0$ if $\dim M=3$. By Theorem \ref{Cor.4.6},
the following definition makes sense:

\begin{enumerate}
\item[(1)]  If $\dim M=3$, the point $p\in M$ is called umbilic whenever 
\[
F_{42}(z,\overline{z},0)=F_{24}(z,\overline{z},0)=0.
\]

\item[(2)]  If $\dim M\geq 5$, the point $p\in M$ is called umbilic whenever 
\[
F_{22}(z,\overline{z},0)=0.
\]
\end{enumerate}

Let $N_{\sigma }$ be a normalization of $M$ and $\phi _{\sigma }=\varphi
\circ \psi $ be an automorphism of a real hyperquadric. We have a
decomposition of $N_{\sigma }$(cf. \cite{Pa1}): 
\[
N_{\sigma }=\phi \circ E\circ \psi 
\]
where $E$ is a normalization of $\psi \left( M\right) $ with identity
initial value. Then we easily verify for each $k\geq 3$%
\[
N_{\sigma }=\phi _{\sigma }+O_{\times }(k+1)\quad \mathrm{if}\text{ }\mathrm{%
and}\text{ }\mathrm{only}\text{ }\mathrm{if\quad }E=id+O_{\times }(k+1). 
\]

\begin{theorem}
\label{Thm.5.14}Let $M$ be an analytic real hypersurface of dimension $3$ in
normal form. Let $N_{\sigma }$ be a normalization of $M$ such that $a\neq 0$
in $\sigma =(C,a,\rho ,r)$. Then the nonumbilicity of the origin $0\in M$ is
equivalent to the following condition: 
\[
N_{\sigma }=\phi _{\sigma }+O_{\times }(7)\text{\textit{\quad }and\textit{%
\quad }}N_{\sigma }\neq \phi _{\sigma }+O_{\times }(8).
\]
\end{theorem}

\proof
It may suffice to show that the normalization $E$ in $N_{\sigma }=\varphi
\circ E\circ \psi $ satisfies 
\begin{equation}
E=id+O_{\times }(7)\text{\textit{\quad }and\textit{\quad }}E\neq
id+O_{\times }(8)  \tag*{(2.1)}  \label{E}
\end{equation}
if and only if the origin is nonumbilic.

Suppose that $M$ is defined by 
\[
v=z\bar{z}+bz^{4}\bar{z}\,^{2}+\bar{b}z^{2}\bar{z}\,^{4}+cz^{5}\bar{z}^{2}+%
\bar{c}z^{2}\bar{z}^{5}+dz^{4}\bar{z}^{3}+\bar{d}z^{3}\bar{z}^{4}+O(8), 
\]
where 
\[
b,c,d\in \Bbb{C}\text{.} 
\]
By the mapping $\psi $, $M$ is transformed up to weight 7 to a real
hypersurface as follows: 
\begin{eqnarray*}
v &=&z\bar{z}+bz^{4}\bar{z}^{2}+\bar{b}z^{2}\bar{z}^{4}+(c+4ib\bar{a})z^{5}%
\bar{z}^{2}+(\overline{c+4ib\bar{a}})z^{2}\bar{z}^{5} \\
&&+(d+2iba)z^{4}\bar{z}^{3}+(\overline{d+2iba})z^{3}\bar{z}^{4} \\
&&+4bauz^{3}\bar{z}^{2}+4\bar{b}\bar{a}uz^{2}\bar{z}^{3}+2b\bar{a}uz^{4}\bar{%
z}+2\bar{b}auz\bar{z}^{4}+O(8)
\end{eqnarray*}
By Theorem \ref{Cor.3.2}, we obtain the normalization $E$ up to weight 7 as
follows: 
\begin{eqnarray*}
z^{*} &=&z+2\bar{a}bz^{4}w-2iabz^{2}w^{2}-\frac{1}{3i}\bar{a}\bar{b}%
w^{3}+O(7) \\
w^{*} &=&w+\frac{2}{3}abzw^{3}+O(8),
\end{eqnarray*}
and $M$ is transformed up to weight 7 to a real hypersurface as follows: 
\begin{align}
v& =z\bar{z}+bz^{4}\bar{z}^{2}+\bar{b}z^{2}\bar{z}^{4}+(c+2ib\bar{a})z^{5}%
\bar{z}^{2}+(\overline{c+2ib\bar{a}})z^{2}\bar{z}^{5}  \nonumber \\
& +(d+\frac{2}{3}iba)z^{4}\bar{z}^{3}+(\overline{d+\frac{2}{3}iba})z^{3}\bar{%
z}^{4}+O(8)  \tag*{(2.2)}  \label{5.2}
\end{align}
Since $a\neq 0,$ the normalization $E$ satisfies the condition \ref{E} if
and only if $b\neq 0$. This completes the proof.\endproof

\begin{corollary}[Moser]
\label{Cor.5.15}Let $M$ be a real hypersurface of dimension $3$ with a
nonumbilic point $p\in M$. Then there is a normal coordinate with center at $%
p\in M$ such that 
\[
v=z\bar{z}+F_{42}(z,\bar{z},u)+F_{24}(z,\bar{z},u)+\sum_{\min (s,t)\geq
2,s+t\geq 7}F_{st}(z,\bar{z},u)
\]
where 
\begin{equation}
F_{24}(z,\bar{z},0)=z^{2}\overline{z}^{4},\quad \Re \left\{ z^{2}\left( 
\frac{\partial F_{24}}{\partial u}\right) (z,\bar{z},0)\right\} =0,\text{%
\quad }F_{43}(z,\bar{z},0)=0.  \tag*{(2.3)}  \label{5.8}
\end{equation}
Further, if two real hypersurfaces are of the reduced normal form and they
are biholomorphic near the origin, then they are related by a mapping as
follows: 
\[
z^{*}=\pm z,\quad w^{*}=w.
\]
The same is true when we replace the condition \ref{5.8} by the following
condition: 
\[
F_{24}(z,\bar{z},0)=z^{2}\overline{z}^{4},\text{\quad }\Re \left\{
z^{2}\left( \frac{\partial F_{24}}{\partial u}\right) (z,\bar{z},0)\right\}
=0,\text{\quad }F_{52}(z,\bar{z},0)=0.
\]
\end{corollary}

\proof
Since the point $p$ is nonumbilic, we take a normalization $N_{\sigma }$(cf.
Theorem \ref{Thm.5.14}) with initial value 
\begin{gather*}
\sigma =(1,a,1,0) \\
a=\frac{3id}{2b}\quad \left( \text{resp.\textit{\quad }}a=-\frac{i\overline{c%
}}{2\overline{b}}\right)
\end{gather*}
where we assume that 
\begin{eqnarray*}
\left. F_{42}\right| _{u=0} &=&bz^{4}\overline{z}^{2}\neq 0 \\
\left. F_{52}\right| _{u=0} &=&cz^{5}\overline{z}^{2} \\
\left. F_{43}\right| _{u=0} &=&dz^{4}\overline{z}^{3}.
\end{eqnarray*}
Then we obtain 
\begin{equation}
F_{43}^{*}\mid _{u=0}=0\quad \left( \text{resp.\textit{\quad }}%
F_{52}^{*}\mid _{u=0}=0\right) .  \tag*{(2.4)}  \label{F43}
\end{equation}
Note that $a=0$ necessarily if $F_{43}\mid _{u=0}=0\left( \text{resp.\textit{%
\quad }}F_{52}\mid _{u=0}=0\right) .$

Suppose that $F_{43}\mid _{u=0}=0\left( \text{resp.\textit{\quad }}%
F_{52}\mid _{u=0}=0\right) .$ Then we take a normalization $N_{\sigma }$
with initial value $\sigma =(\alpha ,0,\alpha \overline{\alpha },0)$: 
\[
z^{*}=\alpha z,\ w^{*}=\alpha \overline{\alpha }w, 
\]
where 
\[
\alpha =\pm \left( \frac{b^{6}}{\bar{b}^{2}}\right) ^{1/4}. 
\]
Then we obtain 
\[
F_{24}^{*}\mid _{u=0}=z^{2}\overline{z}^{4}. 
\]
Note that $a=0$ and $\alpha =\pm 1$ necessarily if 
\begin{gather*}
F_{43}\mid _{u=0}=0\text{ and }F_{24}\mid _{u=0}=z^{2}\overline{z}^{4} \\
\left( \text{resp.\textit{\quad }}F_{52}\mid _{u=0}=0\text{ and }F_{24}\mid
_{u=0}=z^{2}\overline{z}^{4}\right) .
\end{gather*}
Then we carry out another normalization $N_{\sigma }$ with initial value $%
\sigma =(1,0,1,r)$: 
\[
z^{*}=\frac{z}{1-rw},\ \ \ \ \ w^{*}=\frac{w}{1-rw}, 
\]
on a real hypersurface $M$ defined up to weight $8$ by the following
equation: 
\[
v=z\bar{z}+z^{4}\bar{z}^{2}+z^{2}\bar{z}^{4}+duz^{4}\bar{z}^{2}+\bar{d}uz^{2}%
\bar{z}^{4}+O(\left| z\right| ^{7})+O(9). 
\]
Then $M$ is transformed up to weight $8$ to a real hypersurface as follows: 
\[
v=z\bar{z}+z^{4}\bar{z}^{2}+z^{2}\bar{z}^{4}+(d-4r)uz^{4}\bar{z}^{2}+(\bar{d}%
-4r)uz^{2}\bar{z}^{4}+O(\left| z\right| ^{7})+O(9). 
\]
Taking 
\[
r=\Re \left( \frac{d}{4}\right) 
\]
yields 
\[
\Re \left( \left. \frac{\partial F_{24}^{*}}{\partial u}\right|
_{u=0}\right) =0. 
\]

Thus there is a normal coordinate at a nonumbilic point on $M$ such that 
\begin{gather*}
\left. F_{24}\right| _{u=0}=z^{2}\overline{z}^{4},\quad \Re \left( \left.
z^{2}\frac{\partial F_{24}}{\partial u}\right| _{u=0}\right) =0,\text{\quad }%
\left. F_{43}\right| _{u=0}=0 \\
\left( \text{resp.}\quad \left. F_{24}\right| _{u=0}=z^{2}\overline{z}^{4},%
\text{\quad }\Re \left( \left. z^{2}\frac{\partial F_{24}}{\partial u}%
\right| _{u=0}\right) =0,\text{\quad }\left. F_{52}\right| _{u=0}=0.\right)
\end{gather*}
Clearly, a normalization $N_{\sigma }$ between these normal coordinates has
its initial value such as 
\[
\sigma =(\pm 1,0,1,0). 
\]
This completes the proof.\endproof

\textbf{II}. We keep the following conventions unless explicitly specified
otherwise.

\begin{enumerate}
\item[(1)]  All Greek indices run from $1$ to $n$.

\item[(2)]  The summation convention over any repeated Greek indices.

\item[(3)]  Complex conjugation may be indicated by bars over Greek indices.
\ For instance, $z^{\bar{\alpha}}\equiv \overline{z^{\alpha }},$ $c_{\bar{%
\alpha}\bar{\beta}}\equiv \overline{c_{\alpha \beta }}.$ We may use $%
\overline{z}^{\alpha }$ for $\overline{z^{\alpha }}$ with a Greek index, but
may not with a specific numerical value as in $\left( \overline{z}^{\alpha
}\right) =\left( \overline{z^{1}},\cdots ,\overline{z^{n}}\right) .$ We
reserve $\overline{z}^{k}$ for $(\overline{z})^{k},$ $k\in \Bbb{N},$ or for
a notational abbreviation such as $O(z^{s}\overline{z}^{t}),$ $s,t\in \Bbb{N}%
.$

\item[(4)]  We shall raise and lower Greek indices by using the following
matrix 
\[
\text{diag}\{\underbrace{1,\cdots ,1}_{e},\underbrace{-1,\cdots ,-1}_{n-e}\}
\]
where $(e,n-e)$ is the signature of the quadric $\langle z,z\rangle .$ If
necessary, we shall indicate the locations of upper indices by dots in order
as in $A_{\alpha \beta ..}^{\text{ }\gamma \delta }.$
\end{enumerate}

We shall examine the normalization $E$ in the decomposition $N_{\sigma
}=\varphi \circ E\circ \psi $ in low order terms for the case of $\dim M\geq
5.$

\begin{lemma}
\label{compute}Let $M$ be an analytic real hypersurface of dimension$\geq 5$
in normal form, which is defined up to weight $6$ by the following equation: 
\begin{eqnarray*}
v=\langle z,z\rangle +A_{\alpha \beta \bar{\gamma}\bar{\delta}}z^{\alpha
}z^{\beta }z^{\bar{\gamma}}z^{\bar{\delta}}+uB_{\alpha \beta \overline{%
\gamma }\bar{\delta}}z^{\alpha }z^{\beta }z^{\overline{\gamma }}z^{\bar{%
\delta}} \\
+C_{\alpha \beta \gamma \bar{\delta}\bar{\eta}}z^{\alpha }z^{\beta
}z^{\gamma }z^{\bar{\delta}}z^{\bar{\eta}}+C_{\alpha \beta \bar{\gamma}%
\overline{\delta }\overline{\eta }}z^{\alpha }z^{\beta }z^{\bar{\gamma}}z^{%
\overline{\delta }}z^{\overline{\eta }} \\
+D_{\alpha \beta \gamma \delta \bar{\eta}\overline{\xi }}z^{\alpha }z^{\beta
}z^{\gamma }z^{\delta }z^{\bar{\eta}}z^{\overline{\xi }}+O(z^{2}\overline{z}%
^{4})+O(z^{4}\overline{z}^{2})+O(6),
\end{eqnarray*}
where 
\begin{gather*}
\overline{A_{\alpha \beta \bar{\gamma}\bar{\delta}}}=A_{\gamma \delta \bar{%
\alpha}\bar{\beta}},\text{ }\overline{B_{\alpha \beta \overline{\gamma }\bar{%
\delta}}}=B_{\gamma \delta \bar{\alpha}\bar{\beta}},\text{ }\overline{%
C_{\alpha \beta \gamma \bar{\delta}\overline{\eta }}}=C_{\delta \eta \bar{%
\alpha}\bar{\beta}\overline{\gamma }},\text{ }\overline{D_{\alpha \beta
\gamma \overline{\delta }\overline{\eta }\overline{\xi }}}=D_{\delta \eta
\xi \bar{\alpha}\bar{\beta}\overline{\gamma }}, \\
A_{\alpha \beta .\bar{\delta}}^{\text{ }\alpha }=B_{\alpha \beta .\bar{\delta%
}}^{\text{ }\alpha }=0,\text{ }C_{\alpha \beta ..\overline{\eta }}^{\text{ }%
\alpha \beta }=C_{\alpha \beta \gamma ..}^{\text{ }\alpha \beta }=0,\text{ }%
D_{\alpha \beta \gamma ...}^{\text{ }\alpha \beta \gamma }=0,
\end{gather*}
and all barred and unbarred indices are respectively symmetric. Let $%
N_{\sigma }=\varphi \circ E\circ \psi $ be a normalization with initial
value $\sigma =(C,a,\rho ,r)$. Then the normalizing mapping $E$ is given up
to weight $6$ as follows: 
\begin{eqnarray*}
z^{*\alpha }=z^{\alpha }+2g^{\alpha \bar{\delta}}A_{\beta \gamma \bar{\delta}%
\bar{\eta}}a^{\bar{\eta}}z^{\beta }z^{\gamma }w+4ig^{\alpha \bar{\zeta}%
}A_{\beta \gamma \bar{\eta}\bar{\zeta}}a^{\bar{\eta}}z^{\beta }z^{\gamma
}\langle z,a\rangle w \\
+2g^{\alpha \overline{\eta }}C_{\beta \gamma \delta \bar{\eta}\overline{\xi }%
}z^{\beta }z^{\gamma }z^{\delta }a^{\overline{\xi }}w-8iz^{\alpha }A_{\beta
\gamma \bar{\eta}\overline{\xi }}z^{\beta }z^{\gamma }a^{\bar{\eta}}a^{%
\overline{\xi }}w \\
+2g^{\alpha \bar{\delta}}A_{\beta \gamma \bar{\delta}\bar{\eta}}z^{\beta
}a^{\gamma }a^{\overline{\eta }}w^{2} \\
-\frac{3i}{n+2}\cdot g^{\alpha \bar{\delta}}\{C_{\eta \beta \gamma .%
\overline{\delta }}^{\text{ }\eta }z^{\beta }a^{\gamma }+C_{\eta \beta .%
\overline{\delta }\overline{\gamma }}^{\text{ }\eta }z^{\beta }a^{\overline{%
\gamma }}\}w^{2}+O(6) \\
w^{*}=w-4iA_{\alpha \beta \overline{\gamma }\bar{\delta}}z^{\alpha }z^{\beta
}a^{\overline{\gamma }}a^{\bar{\delta}}w^{2}+O(7).
\end{eqnarray*}
\end{lemma}

\proof
By the mapping $\psi $ in the decomposition of $N_{\sigma }=\varphi \circ
E\circ \psi ,$ $M$ is transformed up to weight $6$ to a real hypersurface as
follows: 
\begin{eqnarray*}
v &=&\langle z,z\rangle +F_{02}(z,\overline{z},u)+F_{20}(z,\overline{z},u) \\
&&+F_{12}(z,\overline{z},u)+F_{21}(z,\overline{z},u)+F_{13}(z,\overline{z}%
,u)+F_{31}(z,\overline{z},u) \\
&&+F_{22}(z,\overline{z},u)+F_{23}(z,\overline{z},u)+F_{32}(z,\overline{z},u)
\\
&&+F_{33}(z,\overline{z},u)+O(z^{2}\overline{z}^{4})+O(z^{4}\overline{z}%
^{2})+O(7),
\end{eqnarray*}
where 
\begin{eqnarray*}
F_{02}(z,\overline{z},u) &=&4u^{2}A_{\alpha \beta \overline{\gamma }\bar{%
\delta}}a^{\alpha }a^{\beta }z^{\overline{\gamma }}z^{\bar{\delta}} \\
F_{12}(z,\overline{z},u) &=&2uA_{\alpha \beta \bar{\gamma}\bar{\delta}%
}z^{\alpha }a^{\beta }z^{\bar{\gamma}}z^{\bar{\delta}} \\
F_{13}(z,\overline{z},u) &=&-4iu\langle a,z\rangle A_{\alpha \beta \overline{%
\gamma }\bar{\delta}}z^{\alpha }a^{\beta }z^{\overline{\gamma }}z^{\bar{%
\delta}}+4iu\langle z,z\rangle A_{\alpha \beta \overline{\gamma }\bar{\delta}%
}a^{\alpha }a^{\beta }z^{\overline{\gamma }}z^{\bar{\delta}} \\
&&+2uC_{\alpha \beta \bar{\gamma}\bar{\delta}\overline{\eta }}z^{\alpha
}a^{\beta }z^{\bar{\gamma}}z^{\bar{\delta}}z^{\overline{\eta }} \\
F_{22}(z,\overline{z},u) &=&A_{\alpha \beta \overline{\gamma }\bar{\delta}%
}z^{\alpha }z^{\beta }z^{\overline{\gamma }}z^{\bar{\delta}}+4iu\langle
z,a\rangle A_{\alpha \beta \overline{\gamma }\bar{\delta}}z^{\alpha
}a^{\beta }z^{\overline{\gamma }}z^{\bar{\delta}} \\
&&-4iu\langle a,z\rangle A_{\alpha \beta \overline{\gamma }\bar{\delta}%
}z^{\alpha }z^{\beta }z^{\overline{\gamma }}a^{\bar{\delta}}+uB_{\alpha
\beta \overline{\gamma }\bar{\delta}}z^{\alpha }z^{\beta }z^{\overline{%
\gamma }}z^{\bar{\delta}} \\
&&+3uC_{\alpha \beta \gamma \bar{\delta}\overline{\eta }}z^{\alpha }z^{\beta
}a^{\gamma }z^{\bar{\delta}}z^{\overline{\eta }}+3uC_{\alpha \beta \overline{%
\gamma }\bar{\delta}\overline{\eta }}z^{\alpha }z^{\beta }z^{\overline{%
\gamma }}z^{\bar{\delta}}a^{\overline{\eta }} \\
F_{23}(z,\overline{z},u) &=&C_{\alpha \beta \overline{\gamma }\bar{\delta}%
\overline{\eta }}z^{\alpha }z^{\beta }z^{\overline{\gamma }}z^{\bar{\delta}%
}z^{\overline{\eta }}-2i\langle a,z\rangle A_{\alpha \beta \overline{\gamma }%
\bar{\delta}}z^{\alpha }z^{\beta }z^{\overline{\gamma }}z^{\bar{\delta}} \\
&&+2i\langle z,z\rangle A_{\alpha \beta \overline{\gamma }\bar{\delta}%
}z^{\alpha }a^{\beta }z^{\overline{\gamma }}z^{\bar{\delta}} \\
F_{33}(z,\overline{z},u) &=&D_{\alpha \beta \gamma \bar{\delta}\overline{%
\eta }\overline{\xi }}z^{\alpha }z^{\beta }z^{\gamma }z^{\bar{\delta}}z^{%
\overline{\eta }}z^{\overline{\xi }}-2\langle a,a\rangle \langle z,z\rangle
A_{\alpha \beta \overline{\gamma }\bar{\delta}}z^{\alpha }z^{\beta }z^{%
\overline{\gamma }}z^{\bar{\delta}} \\
&&+4\langle z,a\rangle \langle a,z\rangle A_{\alpha \beta \overline{\gamma }%
\bar{\delta}}z^{\alpha }z^{\beta }z^{\overline{\gamma }}z^{\bar{\delta}%
}-4\langle z,z\rangle \langle z,a\rangle A_{\alpha \beta \overline{\gamma }%
\bar{\delta}}z^{\alpha }a^{\beta }z^{\overline{\gamma }}z^{\bar{\delta}} \\
&&-4\langle z,z\rangle \langle a,z\rangle A_{\alpha \beta \overline{\gamma }%
\bar{\delta}}z^{\alpha }z^{\beta }z^{\overline{\gamma }}a^{\bar{\delta}%
}+4\langle z,z\rangle ^{2}A_{\alpha \beta \overline{\gamma }\bar{\delta}%
}z^{\alpha }a^{\beta }z^{\overline{\gamma }}a^{\bar{\delta}} \\
&&-2i\langle a,z\rangle C_{\alpha \beta \gamma \bar{\delta}\overline{\eta }%
}z^{\alpha }z^{\beta }z^{\gamma }z^{\bar{\delta}}z^{\overline{\eta }%
}+2i\langle z,a\rangle C_{\alpha \beta \overline{\gamma }\bar{\delta}%
\overline{\eta }}z^{\alpha }z^{\beta }z^{\overline{\gamma }}z^{\bar{\delta}%
}z^{\overline{\eta }} \\
&&+3i\langle z,z\rangle C_{\alpha \beta \gamma \bar{\delta}\overline{\eta }%
}z^{\alpha }z^{\beta }a^{\gamma }z^{\bar{\delta}}z^{\overline{\eta }%
}-3i\langle z,z\rangle C_{\alpha \beta \overline{\gamma }\bar{\delta}%
\overline{\eta }}z^{\alpha }z^{\beta }z^{\overline{\gamma }}z^{\bar{\delta}%
}a^{\overline{\eta }}
\end{eqnarray*}

Since $M$ is in normal form, $E=id+O_{\times }(5)$ by Lemma \ref{Lem.4.1}.
So the function $p(u)$ satisfies 
\[
p(u)=\frac{1}{2}p^{\prime \prime }(0)u^{2}+O(6). 
\]
Let $g(z,w)$ be a holomorphic function(cf. Lemma 3 in the paper \cite{Pa1})
defined by 
\begin{eqnarray*}
g(z,w)-g(0,w) &=&-2iF(p(w),\overline{p}(w),w) \\
&&+2iF\left( z+p(w),\overline{p}(w),w+\frac{1}{2}\{g(z,w)-g(0,w)\}\right) \\
g(0,w) &=&iF(p(w),\overline{p}(w),w).
\end{eqnarray*}
Then the function $g(0,u)$ satisfies 
\[
g(0,u)=iF(p(u),\overline{p}(u),u)=O(8). 
\]
Hence the holomorphic function $g(z,w)$ is implicitly defined up to weight $%
8 $ as follows: 
\[
g=2iF(z+p(w),\overline{p}(w),w+\frac{1}{2}g)+O(8). 
\]
Thus we obtain 
\begin{eqnarray*}
g &=&2iF(z+p(w),\overline{p}(w),w+\frac{1}{2}g)+O(8) \\
&=&i\langle z,p^{\prime \prime }(0)\rangle w^{2}+2iF(z,0,w+\frac{1}{2}g)+O(7)
\\
&=&i\langle z,p^{\prime \prime }(0)\rangle w^{2}+4i(w+\frac{1}{2}%
g)^{2}A_{\alpha \beta \overline{\gamma }\bar{\delta}}z^{\alpha }z^{\beta }a^{%
\overline{\gamma }}a^{\bar{\delta}}+O(7) \\
&=&i\langle z,p^{\prime \prime }(0)\rangle w^{2}+4iw^{2}A_{\alpha \beta 
\overline{\gamma }\bar{\delta}}z^{\alpha }z^{\beta }a^{\overline{\gamma }}a^{%
\bar{\delta}}+O(7).
\end{eqnarray*}

We carry out the following mapping: 
\begin{eqnarray*}
z &=&z^{*}+\frac{1}{2}p^{\prime \prime }(0)w^{*2}+O(6) \\
w &=&w^{*}+i\langle z^{*},p^{\prime \prime }(0)\rangle
w^{*2}+4iw^{*2}A_{\alpha \beta \overline{\gamma }\bar{\delta}}z^{*\alpha
}z^{*\beta }a^{*\overline{\gamma }}a^{*\bar{\delta}}+O(7).
\end{eqnarray*}
Its inverse mapping is given as follows: 
\begin{align}
z^{*}=& z-\frac{1}{2}p^{\prime \prime }(0)w^{2}+O(6)  \nonumber \\
w^{*}=& w-i\langle z,p^{\prime \prime }(0)\rangle w^{2}-4iw^{2}A_{\alpha
\beta \overline{\gamma }\bar{\delta}}z^{\alpha }z^{\beta }a^{\overline{%
\gamma }}a^{\bar{\delta}}+O(7).  \tag*{(2.5)}  \label{pope}
\end{align}
Then $M$ is transformed by the mapping \ref{pope} to a real hypersurface up
to weight $6$ as follows: 
\begin{eqnarray*}
v &=&\langle z,z\rangle -2iu\langle z,p^{\prime \prime }(0)\rangle \langle
z,z\rangle +2iu\langle p^{\prime \prime }(0),z\rangle \langle z,z\rangle \\
&&-4iu\langle z,z\rangle A_{\alpha \beta \overline{\gamma }\bar{\delta}%
}z^{\alpha }z^{\beta }a^{\overline{\gamma }}a^{\bar{\delta}}+4iu\langle
z,z\rangle A_{\alpha \beta \overline{\gamma }\bar{\delta}}a^{\alpha
}a^{\beta }z^{\overline{\gamma }}z^{\bar{\delta}} \\
&&+F_{12}(z,\overline{z},u)+F_{21}(z,\overline{z},u)+F_{13}(z,\overline{z}%
,u)+F_{31}(z,\overline{z},u) \\
&&+F_{22}(z,\overline{z},u)+F_{23}(z,\overline{z},u)+F_{32}(z,\overline{z},u)
\\
&&+F_{33}(z,\overline{z},u)+O(z^{2}\overline{z}^{4})+O(z^{4}\overline{z}%
^{2})+O(7) \\
&=&\langle z,z\rangle +F_{12}^{*}(z,\overline{z},u)+F_{21}^{*}(z,\overline{z}%
,u)+F_{13}(z,\overline{z},u)+F_{31}(z,\overline{z},u) \\
&&+F_{22}(z,\overline{z},u)+F_{23}^{*}(z,\overline{z},u)+F_{32}^{*}(z,%
\overline{z},u) \\
&&+F_{33}(z,\overline{z},u)+O(z^{2}\overline{z}^{4})+O(z^{4}\overline{z}%
^{2})+O(7),
\end{eqnarray*}
where 
\begin{eqnarray*}
F_{12}^{*}(z,\overline{z},u) &=&F_{12}(z,\overline{z},u)+2iu\langle
p^{\prime \prime }(0),z\rangle \langle z,z\rangle \\
F_{13}^{*}(z,\overline{z},u) &=&F_{13}(z,\overline{z},u)+4i\langle
z,z\rangle A_{\alpha \beta \overline{\gamma }\bar{\delta}}a^{\alpha
}a^{\beta }z^{\overline{\gamma }}z^{\bar{\delta}}.
\end{eqnarray*}

We carry out the following mapping: 
\begin{align}
z^{*\alpha }& =z^{\alpha }-2iz^{\alpha }\langle z,p^{\prime \prime
}(0)\rangle w+2g^{\alpha \bar{\delta}}A_{\beta \gamma \bar{\delta}\bar{\eta}%
}a^{\bar{\eta}}z^{\beta }z^{\gamma }w  \nonumber \\
& +4ig^{\alpha \bar{\zeta}}A_{\beta \gamma \bar{\eta}\bar{\zeta}}a^{\bar{\eta%
}}z^{\beta }z^{\gamma }\langle z,a\rangle w+2g^{\alpha \overline{\eta }%
}C_{\beta \gamma \delta \bar{\eta}\overline{\xi }}z^{\beta }z^{\gamma
}z^{\delta }a^{\overline{\xi }}w  \nonumber \\
& -8iz^{\alpha }A_{\beta \gamma \bar{\eta}\overline{\xi }}z^{\beta
}z^{\gamma }a^{\bar{\eta}}a^{\overline{\xi }}w+O(6)  \nonumber \\
w^{*}& =w.  \tag*{(2.6)}  \label{5.3}
\end{align}
Then the real hypersurface is transformed to 
\begin{eqnarray*}
v &=&\langle z,z\rangle +4u^{2}A_{\alpha \beta \overline{\gamma }\bar{\delta}%
}z^{\alpha }a^{\beta }z^{\overline{\gamma }}a^{\bar{\delta}} \\
&&+F_{22}(z,\overline{z},u)+F_{23}^{*}(z,\overline{z},u)+F_{32}^{*}(z,%
\overline{z},u)+F_{33}(z,\overline{z},u) \\
&&+O(z^{2}\overline{z}^{4})+O(z^{4}\overline{z}^{2})+O(7),
\end{eqnarray*}
where 
\[
F_{23}^{*}(z,\overline{z},u)=F_{23}(z,\overline{z},u)-2\langle z,z\rangle
^{2}\langle p^{\prime \prime }(0),z\rangle +2i\langle z,z\rangle A_{\alpha
\beta \overline{\gamma }\bar{\delta}}z^{\alpha }a^{\beta }z^{\overline{%
\gamma }}z^{\bar{\delta}}. 
\]
The condition $\Delta ^{2}F_{23}^{*}=0$ determines the coefficient $%
p^{\prime \prime }(0)$ so that 
\[
p^{\prime \prime }(0)=0. 
\]

Then we take a matrix $L$ defined by 
\[
\langle Lz,z\rangle =2A_{\alpha \beta \overline{\gamma }\bar{\delta}%
}z^{\alpha }a^{\beta }z^{\overline{\gamma }}a^{\bar{\delta}}. 
\]
Since $v=\langle z,z\rangle +O(4)$ on the real hypersurface, we have 
\begin{eqnarray*}
\langle z,z\rangle +4u^{2}A_{\alpha \beta \overline{\gamma }\bar{\delta}%
}z^{\alpha }a^{\beta }z^{\overline{\gamma }}a^{\bar{\delta}} &=&\langle
z+u^{2}Lz,z+u^{2}Lz\rangle +O(7) \\
&=&\langle z+w^{2}Lz,z+w^{2}Lz\rangle \\
&&+4\langle z,z\rangle ^{2}A_{\alpha \beta \overline{\gamma }\bar{\delta}%
}z^{\alpha }a^{\beta }z^{\overline{\gamma }}a^{\bar{\delta}}+O(7).
\end{eqnarray*}
We consider matrices $A,B$ such that 
\[
\langle z,z\rangle =\langle z+uAz+u^{2}Bz,z+uAz+u^{2}Bz\rangle +O(7). 
\]
Since $E=id+O_{\times }(5),$ we have $A=0$ so that 
\begin{eqnarray*}
\langle z,z\rangle &=&\langle z+u^{2}Bz,z+u^{2}Bz\rangle +O(7) \\
&=&\langle z+w^{2}Bz,z+w^{2}Bz\rangle -4iu\langle z,z\rangle \langle
Bz,z\rangle +O(7),
\end{eqnarray*}
on the real hypersurface, where 
\[
\langle Bz,z\rangle +\langle z,Bz\rangle =0. 
\]

We carry out the mapping: 
\begin{align}
z^{*}=& z+Lzw^{2}+Bzw^{2}  \nonumber \\
w^{*}=& w,  \tag*{(2.7)}  \label{third}
\end{align}
Then the real hypersurface $M$ is transformed up to weight $6$ to 
\begin{eqnarray*}
v &=&\langle z,z\rangle +F_{22}^{*}(z,\overline{z},u)+F_{23}^{*}(z,\overline{%
z},u)+F_{32}^{*}(z,\overline{z},u)+F_{33}^{*}(z,\overline{z},u) \\
&&+O(z^{2}\overline{z}^{4})+O(z^{4}\overline{z}^{2})+O(7),
\end{eqnarray*}
where 
\begin{eqnarray*}
F_{22}^{*}(z,\overline{z},u) &=&F_{22}(z,\overline{z},u)-4iu\langle
z,z\rangle \langle Bz,z\rangle \\
F_{33}^{*}(z,\overline{z},u) &=&F_{33}(z,\overline{z},u)+4\langle z,z\rangle
^{2}A_{\alpha \beta \overline{\gamma }\bar{\delta}}z^{\alpha }a^{\beta }z^{%
\overline{\gamma }}a^{\bar{\delta}}.
\end{eqnarray*}
The condition $\Delta ^{2}F_{22}=0$ determines the matrix $B$ so that 
\[
\langle Bz,z\rangle =-\frac{3i}{n+2}\{C_{\delta \alpha \beta .\overline{%
\gamma }}^{\text{ }\delta }z^{\alpha }a^{\beta }z^{\overline{\gamma }%
}+C_{\delta \alpha .\overline{\beta }\overline{\gamma }}^{\text{ }\delta
}z^{\alpha }z^{\overline{\beta }}a^{\overline{\gamma }}\}. 
\]
We easily verify $\Delta ^{3}F_{33}^{*}=0$ by using $A_{\alpha \beta .%
\overline{\delta }}^{\text{ }\alpha }=0$.

Thus we have showed by composing the mappings \ref{pope}, \ref{5.3}, \ref
{third} that the normalization $E$ is given up to weight $6$ by the
following mapping: 
\begin{eqnarray*}
z^{*\alpha } &=&z^{\alpha }+2g^{\alpha \bar{\delta}}A_{\beta \gamma \bar{%
\delta}\bar{\eta}}a^{\bar{\eta}}z^{\beta }z^{\gamma }w+4ig^{\alpha \bar{\zeta%
}}A_{\beta \gamma \bar{\eta}\bar{\zeta}}a^{\bar{\eta}}z^{\beta }z^{\gamma
}\langle z,a\rangle w \\
&&+2g^{\alpha \overline{\eta }}C_{\beta \gamma \delta \bar{\eta}\overline{%
\xi }}z^{\beta }z^{\gamma }z^{\delta }a^{\overline{\xi }}w-8iz^{\alpha
}A_{\beta \gamma \bar{\eta}\overline{\xi }}z^{\beta }z^{\gamma }a^{\bar{\eta}%
}a^{\overline{\xi }}w \\
&&+2g^{\alpha \bar{\delta}}A_{\beta \gamma \bar{\delta}\bar{\eta}}z^{\beta
}a^{\gamma }a^{\overline{\eta }}w^{2} \\
&&-\frac{3i}{n+2}\cdot g^{\alpha \bar{\delta}}\{C_{\eta \beta \gamma .%
\overline{\delta }}^{\text{ }\eta }z^{\beta }a^{\gamma }+C_{\eta \beta .%
\overline{\delta }\overline{\gamma }}^{\text{ }\eta }z^{\beta }a^{\overline{%
\gamma }}\}w^{2}+O(6) \\
w^{*} &=&w-4iA_{\alpha \beta \overline{\gamma }\bar{\delta}}z^{\alpha
}z^{\beta }a^{\overline{\gamma }}a^{\bar{\delta}}w^{2}+O(7).
\end{eqnarray*}
This completes the proof.\endproof

\begin{theorem}
\label{Th5}Let $M$ be a real hypersurface of dimension$\geq 5$ in normal
form defined by 
\[
v=\langle z,z\rangle +\sum_{\min (s,t)\geq 2}F_{st}(z,\bar{z},u).
\]
Suppose that the origin $0\in M$ is nonumbilic. Let $N_{\sigma }$ be a
normalization such that the parameter $a$ in $\sigma =(C,a,\rho ,r)$
satisfies the condition: 
\begin{equation}
\sum_{\alpha }a^{\alpha }\left( \frac{\partial F_{22}}{\partial z^{\alpha }}%
\right) (z,\overline{z},0)\neq 0.  \tag*{(2.8)}  \label{5.9}
\end{equation}
Then the normalization $N_{\sigma }$ satisfies the following condition: 
\begin{equation}
N_{\sigma }=\phi _{\sigma }+O_{\times }(5)\text{ and }N_{\sigma }\neq \phi
_{\sigma }+O_{\times }(6).  \tag*{(2.9)}  \label{ene}
\end{equation}
Conversely, if there is a normalization $N_{\sigma }$ satisfying the
condition \ref{ene}, then the origin is nonumbilic.
\end{theorem}

\proof
It may suffice to show that the normalization $E$ in $N_{\sigma }=\varphi
\circ E\circ \psi $ satisfies 
\[
E=\phi _{\sigma }+O_{\times }(5)\text{ and }E\neq \phi _{\sigma }+O_{\times
}(6) 
\]
whenever the condition \ref{5.9} is satisfied.

From Lemma \ref{compute}, the normalization $E$ is given up to weight $5$ as
follows: 
\begin{eqnarray*}
z^{*\alpha } &=&z^{\alpha }+2g^{\alpha \bar{\delta}}A_{\beta \gamma \bar{%
\delta}\bar{\eta}}a^{\bar{\eta}}z^{\beta }z^{\gamma }w+O(5) \\
w^{*} &=&w+O(6).
\end{eqnarray*}
Therefore the condition \ref{5.9} implies 
\[
E=\phi _{\sigma }+O_{\times }(5)\text{ and }E\neq \phi _{\sigma }+O_{\times
}(6) 
\]
which is equivalent to 
\[
N_{\sigma }=\phi _{\sigma }+O_{\times }(5)\text{ and }N_{\sigma }\neq \phi
_{\sigma }+O_{\times }(6). 
\]
Clearly, the origin is nonumbilic, i.e., 
\[
A_{\alpha \beta \overline{\gamma }\overline{\delta }}z^{\alpha }z^{\beta }z^{%
\overline{\gamma }}z^{\overline{\delta }}\neq 0 
\]
whenever there is a normalization $N_{\sigma }$ such that 
\[
N_{\sigma }=\phi _{\sigma }+O_{\times }(5)\text{ and }N_{\sigma }\neq \phi
_{\sigma }+O_{\times }(6). 
\]
This completes the proof.\endproof

\begin{corollary}[Webster]
\label{Cor.5.16}Let $M$ be an analytic real hypersurface of dimension$\geq 5$
in a normal coordinate with center at a point $p\in M$ as follows: 
\[
v=\langle z,z\rangle +\sum_{\min (s,t)\geq 2}F_{st}(z,\bar{z},u)
\]
where the function $F_{22}(z,\bar{z},u)$ satisfy the following conditions: 
\[
\left. \Delta ^{4}(F_{22})^{2}\right| _{0}\neq 0.
\]
Then there is a normal coordinate which satisfy the following condition: 
\[
\left. \Delta ^{4}(F_{22})^{2}\right| _{0}=\pm 1,\text{\quad }\left. \frac{d%
}{du}\Delta ^{4}(F_{22})^{2}\right| _{0}=0,\text{\quad }\left. \Delta
^{4}\left( F_{22}\frac{\partial F_{23}}{\partial \overline{z}}\right)
\right| _{0}=0.
\]
Further, if two real hypersurfaces are of the reduced normal form and they
are biholomorphic near the origin, then they are related by a mapping as
follows: 
\[
z^{*}=Cz,\text{\quad }w^{*}=\pm w,
\]
where 
\[
\langle Cz,Cz\rangle =\pm \langle z,z\rangle .
\]
\end{corollary}

\proof
Let $M$ be defined in a normal coordinate up to weight $5$ as follows: 
\begin{eqnarray*}
v &=&\langle z,z\rangle +A_{\alpha \beta \bar{\gamma}\bar{\delta}}z^{\alpha
}z^{\beta }z^{\bar{\gamma}}z^{\bar{\delta}}+C_{\alpha \beta \gamma \bar{%
\delta}\bar{\eta}}z^{\alpha }z^{\beta }z^{\gamma }z^{\bar{\delta}}z^{\bar{%
\eta}} \\
&&+C_{\bar{\alpha}\bar{\beta}\bar{\gamma}\delta \eta }z^{\bar{\alpha}}z^{%
\bar{\beta}}z^{\bar{\gamma}}z^{\delta }z^{\eta }+O(6).
\end{eqnarray*}
By the mapping $E\circ \psi $, $M$ is mapped up to weight $5$ in Theorem \ref
{Th5} to 
\begin{eqnarray*}
v &=&\langle z,z\rangle +A_{\alpha \beta \bar{\gamma}\bar{\delta}%
}^{*}z^{\alpha }z^{\beta }z^{\bar{\gamma}}z^{\bar{\delta}}+C_{\alpha \beta
\gamma \bar{\delta}\bar{\eta}}^{*}z^{\alpha }z^{\beta }z^{\gamma }z^{\bar{%
\delta}}z^{\bar{\eta}} \\
&&+C_{\bar{\alpha}\bar{\beta}\bar{\gamma}\delta \eta }^{*}z^{\bar{\alpha}}z^{%
\bar{\beta}}z^{\bar{\gamma}}z^{\delta }z^{\eta }+O(6),
\end{eqnarray*}
where 
\begin{eqnarray*}
A_{\alpha \beta \bar{\gamma}\bar{\delta}}^{*} &=&A_{\alpha \beta \bar{\gamma}%
\bar{\delta}}, \\
C_{\alpha \beta \gamma \bar{\delta}\bar{\eta}}^{*} &=&C_{\alpha \beta \gamma 
\bar{\delta}\bar{\eta}}-\frac{i}{3}\left\{ g_{\alpha \bar{\delta}}A_{\beta
\gamma \bar{\eta}\overline{\xi }}a^{\overline{\xi }}+g_{\beta \bar{\delta}%
}A_{\gamma \alpha \bar{\eta}\overline{\xi }}a^{\overline{\xi }}+g_{\gamma 
\bar{\delta}}A_{\alpha \beta \bar{\eta}\overline{\xi }}a^{\overline{\xi }%
}\right. \\
&&\hspace{2cm}\left. +g_{\alpha \bar{\eta}}A_{\beta \gamma \bar{\delta}%
\overline{\xi }}a^{\overline{\xi }}+g_{\beta \bar{\eta}}A_{\gamma \alpha 
\bar{\delta}\overline{\xi }}a^{\overline{\xi }}+g_{\gamma \bar{\eta}%
}A_{\alpha \beta \bar{\delta}\overline{\xi }}a^{\overline{\xi }}\right\} .
\end{eqnarray*}
Then we obtain 
\begin{equation}
A_{\alpha \beta ..}^{*\text{ }\gamma \delta }C_{\gamma \delta \zeta ..}^{*%
\text{\ }\alpha \beta }=A_{\alpha \beta ..}^{\text{ }\gamma \delta
}C_{\gamma \delta \zeta ..}^{\text{ }\alpha \beta }-\frac{2i}{3}(A_{\alpha
\beta ..}^{\text{ \ }\gamma \delta }A_{\gamma \delta ..}^{\text{ \ }\alpha
\beta })g_{\zeta \bar{\eta}}a^{\bar{\eta}},  \tag*{(2.10)}  \label{a0}
\end{equation}
where we raise indices by using the following matrix: 
\[
(g_{\beta \bar{\delta}})=\text{diag}\{\underbrace{1,\cdots ,1}_{e},%
\underbrace{-1,\cdots ,-1}_{n-e}\}. 
\]
Suppose that 
\begin{eqnarray*}
F_{22}(z,\overline{z},u) &=&N_{\alpha \beta \overline{\gamma }\overline{%
\delta }}z^{\alpha }z^{\beta }z^{\overline{\gamma }}z^{\overline{\delta }} \\
F_{23}(z,\overline{z},u) &=&N_{\alpha \beta \overline{\gamma }\overline{%
\delta }\overline{\eta }}z^{\alpha }z^{\beta }z^{\overline{\gamma }}z^{%
\overline{\delta }}z^{\overline{\eta }}.
\end{eqnarray*}
Since $\Delta F_{22}=0,$ we obtain 
\begin{eqnarray*}
\Delta ^{4}(F_{22})^{2} &=&3\cdot 2^{5}N_{\alpha \beta ..}^{\text{ \ }\gamma
\delta }N_{\gamma \delta ..}^{\text{ \ }\alpha \beta } \\
\Delta ^{4}\left( F_{22}\frac{\partial F_{23}}{\partial \overline{z}^{\zeta }%
}\right) &=&9\cdot 2^{5}N_{\alpha \beta ..}^{\text{ }\gamma \delta
}N_{\gamma \delta ..\overline{\zeta }}^{\text{ }\alpha \beta }.
\end{eqnarray*}
Thus the condition equality \ref{a0} reads 
\[
\left. \Delta ^{4}\left( F_{22}^{*}\frac{\partial F_{23}^{*}}{\partial 
\overline{z}}\right) \right| _{0}=\left. \Delta ^{4}\left( F_{22}\frac{%
\partial F_{23}}{\partial \overline{z}}\right) \right| _{0}+2ia\left. \Delta
^{4}(F_{22})^{2}\right| _{0}. 
\]
We take a normalization $N_{\sigma }$ with 
\begin{gather*}
\sigma =(id_{n\times n},a,1,0) \\
a=\frac{i}{2}\frac{\left. \Delta ^{4}\left( F_{22}\frac{\partial F_{23}}{%
\partial \overline{z}}\right) \right| _{0}}{\left. \Delta
^{4}(F_{22})^{2}\right| _{0}}
\end{gather*}
so that 
\[
\left. \Delta ^{4}\left( F_{22}^{*}\frac{\partial F_{23}^{*}}{\partial 
\overline{z}}\right) \right| _{0}=0. 
\]

Note that $a=0$ necessarily if 
\[
\left. \Delta ^{4}\left( F_{22}\frac{\partial F_{23}}{\partial \overline{z}}%
\right) \right| _{0}=0. 
\]
Then we carry out a normalization $N_{\sigma }$ with $\sigma =(\sqrt{\rho }%
,0,\rho ,r),$ $\rho >0$: 
\[
z^{*}=\frac{\sqrt{\rho }z}{1-rw},\ \ \ \ \ w^{*}=\frac{\rho w}{1-rw} 
\]
so that 
\begin{eqnarray*}
F_{22}^{*}(z,\overline{z},0) &=&\rho ^{-1}F_{22}(z,\overline{z},0) \\
\left( \frac{\partial F_{22}^{*}}{\partial u}\right) (z,\overline{z},0)
&=&\rho ^{-1}\left( \frac{\partial F_{22}}{\partial u}\right) (z,\overline{z}%
,0)-2\rho ^{-2}rF_{22}(z,\overline{z},0).
\end{eqnarray*}
Then we obtain 
\begin{eqnarray*}
\left. \Delta ^{4}(F_{22}^{*})^{2}\right| _{0} &=&\rho ^{-2}\left. \Delta
^{4}(F_{22})^{2}\right| _{0} \\
\left. \frac{d}{du}\Delta ^{4}(F_{22}^{*})^{2}\right| _{0} &=&\rho
^{-2}\left. \frac{d}{du}\Delta ^{4}(F_{22})^{2}\right| _{0}-4\rho
^{-3}r\left. \Delta ^{4}(F_{22})^{2}\right| _{0}.
\end{eqnarray*}
We take 
\begin{eqnarray*}
\rho &=&\left. \sqrt{\left| \Delta ^{4}(F_{22})^{2}\right| }\right| _{0} \\
r &=&\text{sign}\left\{ \left. \Delta ^{4}(F_{22})^{2}\right| _{0}\right\} 
\frac{\left. \frac{d}{du}\Delta ^{4}(F_{22})^{2}\right| _{0}}{4\left( \left. 
\sqrt{\left| \Delta ^{4}(F_{22})^{2}\right| }\right| _{0}\right) ^{3}}
\end{eqnarray*}
so that 
\begin{eqnarray*}
\left. \Delta ^{4}(F_{22}^{*})^{2}\right| _{0} &=&\pm 1 \\
\left. \frac{d}{du}\Delta ^{4}(F_{22}^{*})^{2}\right| _{0} &=&0.
\end{eqnarray*}

Suppose that $M,M^{\prime }$ are in reduced normal form by the conditions 
\[
\left. \Delta ^{4}\left( F_{22}\frac{\partial F_{23}}{\partial \overline{z}}%
\right) \right| _{0}=0,\quad \left. \Delta ^{4}(F_{22})^{2}\right| _{0}=\pm
1,\quad \left. \frac{d}{du}\Delta ^{4}(F_{22})^{2}\right| _{0}=0 
\]
and there is a normalization of $M,$ $N_{\sigma },$ with $\sigma =(C,a,\rho
,r)$ satisfying $M^{\prime }=N_{\sigma }(M).$ Then necessarily we have 
\[
a=0,\quad \rho =\pm 1,\quad r=0 
\]
so that 
\[
N_{\sigma }:\left\{ 
\begin{array}{l}
z^{*}=Cz, \\ 
w^{*}=\pm w,
\end{array}
\right. 
\]
where 
\[
\langle Cz,Cz\rangle =\pm \langle z,z\rangle . 
\]
This completes the proof.\endproof

\section{Spherical analytic real hypersurfaces}

\textbf{I}. We shall study a nondegenerate analytic real hypersurface $M$
with an open subset of umbilic points.

\begin{lemma}
\label{Lem.5.19}Let $k$ be a positive integer$\geq 7.$ Suppose that $M$ is a
real hypersurface in normal form such that 
\[
v=\langle z,z\rangle +\sum_{\min (s,t)\geq 2,s+t=k}F_{st}(z,\bar{z}%
,u)+\sum_{\min (s,t)\geq 2,s+t\geq k+1}F_{st}(z,\bar{z},u),
\]
where 
\[
\sum_{\min (s,t)\geq 2,s+t=k}F_{st}(z,\bar{z},0)\neq 0.
\]
Then there is a vector $a\in \Bbb{C}^{n}$ and a normalization of $M,$ $%
N_{\sigma },$ $\sigma =(id_{n\times n},a,1,0),$ such that 
\[
F^{*}(z,\bar{z},u)=\sum_{\min (s,t)\geq 2}^{s+t=k-1}F_{st}^{*}(z,\bar{z}%
,u)+\sum_{\min (s,t)\geq 2,s+t\geq k}F_{st}^{*}(z,\bar{z},u),
\]
where the real hypersurface $N_{\sigma }\left( M\right) $ is defined by 
\[
v=\langle z,z\rangle +F^{*}(z,\bar{z},u)
\]
and 
\[
\sum_{\min (s,t)\geq 2,s+t=k-1}F_{st}^{*}(z,\bar{z},u)\neq 0.
\]
\end{lemma}

\proof
The real hypersurface $M$ is defined up to weight $k$ as follows: 
\begin{eqnarray*}
v &=&\langle z,z\rangle +\sum_{\min (s,t)\geq 2,s+t=k}F_{st}(z,\bar{z},0) \\
&&+\sum_{\min (s,t)\geq 2,s+t=k+1}O(z^{s}\bar{z}^{t})+O(k+2).
\end{eqnarray*}
By the mapping $\psi $ in the decomposition of $N_{\sigma }=\varphi \circ
E\circ \psi ,$ we obtain 
\begin{eqnarray*}
v &=&\langle z,z\rangle \\
&&+\sum_{\min (s,t)\geq 2,s+t=k}\sum_{\alpha }u\left\{ a^{\alpha }\left( 
\frac{\partial F_{st}}{\partial z^{\alpha }}\right) (z,\bar{z},0)+\overline{a%
}^{\alpha }\left( \frac{\partial F_{st}}{\partial \overline{z}^{\alpha }}%
\right) (z,\bar{z},0)\right\} \\
&&+\sum_{\min (s,t)\geq 2,s+t=k}F_{st}(z,\bar{z},0)+\sum_{\min (s,t)\geq
2,s+t=k+1}O(z^{s}\bar{z}^{t})+O(k+2) \\
&=&\langle z,z\rangle \\
&&+\sum_{\alpha }u\left\{ a^{\alpha }\left( \frac{\partial F_{2,k-2}}{%
\partial z^{\alpha }}\right) (z,\bar{z},0)+\overline{a}^{\alpha }\left( 
\frac{\partial F_{k-2,2}}{\partial \overline{z}^{\alpha }}\right) (z,\bar{z}%
,0)\right\} \\
&&+\sum_{\min (s,t)\geq 2,s+t=k-1}\sum_{\alpha }u\left\{ a^{\alpha }\left( 
\frac{\partial F_{s+1,t}}{\partial z^{\alpha }}\right) (z,\bar{z},0)\right.
\\
&&\hspace{5cm}\left. +\overline{a}^{\alpha }\left( \frac{\partial F_{s,t+1}}{%
\partial \overline{z}^{\alpha }}\right) (z,\bar{z},0)\right\} \\
&&+\sum_{\min (s,t)\geq 2,s+t=k}F_{st}(z,\bar{z},0)+\sum_{\min (s,t)\geq
2,s+t=k+1}O(z^{s}\bar{z}^{t})+O(k+2).
\end{eqnarray*}
Then by the normalization $E,$ we obtain for $k\geq 8$%
\begin{eqnarray*}
v &=&\langle z,z\rangle \\
&&+\sum_{\min (s,t)\geq 2,s+t=k-1}\sum_{\alpha }u\left\{ a^{\alpha }\left( 
\frac{\partial F_{s+1,t}}{\partial z^{\alpha }}\right) (z,\bar{z},0)\right.
\\
&&\hspace{5cm}\left. +\overline{a}^{\alpha }\left( \frac{\partial F_{s,t+1}}{%
\partial \overline{z}^{\alpha }}\right) (z,\bar{z},0)\right\} \\
&&+\sum_{\min (s,t)\geq 2,s+t=k}F_{st}(z,\bar{z},0)+\sum_{\min (s,t)\geq
2,s+t=k+1}O(z^{s}\bar{z}^{t})+O(k+2).
\end{eqnarray*}
For $k=7,$ we take a vector $a=(a^{\alpha })\in \Bbb{C}^{n}$ so that 
\begin{align}
\left( \frac{\partial F_{42}^{*}}{\partial u}\right) (z,\overline{z},0)=&
\sum_{\alpha }\left\{ a^{\alpha }\left( \frac{\partial F_{52}}{\partial
z^{\alpha }}\right) (z,\bar{z},0)+\overline{a}^{\alpha }\left( \frac{%
\partial F_{43}}{\partial \overline{z}^{\alpha }}\right) (z,\bar{z}%
,0)\right\}  \nonumber \\
\neq & 0.  \tag*{(3.1)}  \label{a.3}
\end{align}
By the normalization $E,$ we obtain for $k=7$%
\begin{eqnarray*}
v &=&\langle z,z\rangle \\
&&+\sum_{\min (s,t)\geq 2,s+t=6}\sum_{\alpha }u\left\{ a^{\alpha }\left( 
\frac{\partial F_{s+1,t}}{\partial z^{\alpha }}\right) (z,\bar{z},0)\right.
\\
&&\hspace{5cm}\left. +\overline{a}^{\alpha }\left( \frac{\partial F_{s,t+1}}{%
\partial \overline{z}^{\alpha }}\right) (z,\bar{z},0)\right\} \\
&&+\kappa \langle z,z\rangle ^{3}+\sum_{\min (s,t)\geq 2,s+t=7}F_{st}(z,\bar{%
z},0) \\
&&+\sum_{\min (s,t)\geq 2,s+t=8}O(z^{s}\bar{z}^{t})+O(9)
\end{eqnarray*}
where the constant $\kappa $ is determined so that $\Delta ^{3}F_{33}^{*}=0.$

We easily see that the condition \ref{a.3} remains valid for $k=7$. Thus,
for $k\geq 7,$ we have showed that there is a vector $a\in \Bbb{C}^{n}$ such
that 
\[
\sum_{\min (s,t)\geq 2,s+t=k-1}\left( \frac{\partial F_{st}^{*}}{\partial u}%
\right) (z,\bar{z},0)u\neq 0 
\]
This completes the proof.\endproof

\begin{lemma}
\label{new}Let $M$ be a real hypersurface in normal form, which is defined
up to weight $6$ by the following equation: 
\begin{eqnarray*}
v=\langle z,z\rangle +A_{\alpha \beta \gamma \delta \overline{\eta }%
\overline{\xi }}z^{\alpha }z^{\beta }z^{\gamma }z^{\delta }z^{\overline{\eta 
}}z^{\overline{\xi }}+A_{\alpha \beta \gamma \overline{\delta }\overline{%
\eta }\overline{\xi }}z^{\alpha }z^{\beta }z^{\gamma }z^{\overline{\delta }%
}z^{\overline{\eta }}z^{\overline{\xi }} \\
+A_{\alpha \beta \overline{\gamma }\overline{\delta }\overline{\eta }%
\overline{\xi }}z^{\alpha }z^{\beta }z^{\overline{\gamma }}z^{\overline{%
\delta }}z^{\overline{\eta }}z^{\overline{\xi }}+O(7)
\end{eqnarray*}
where 
\[
\overline{A_{\alpha \beta \gamma \delta \overline{\eta }\overline{\xi }}}%
=A_{\eta \xi \overline{\alpha }\overline{\beta }\overline{\gamma }\overline{%
\delta }},\text{\quad }\overline{A_{\alpha \beta \gamma \overline{\delta }%
\overline{\eta }\overline{\xi }}}=A_{\delta \eta \xi \overline{\alpha }%
\overline{\beta }\overline{\gamma }},\text{\quad }A_{\alpha \beta \gamma
...}^{\text{ }\alpha \beta \gamma }=0.
\]

Let $N_{\sigma }$ be a normalization of $M$ with initial value $\sigma
=(id_{n\times n},a,1,0)$. Suppose that $N_{\sigma }\left( M\right) $ is
defined by the equation 
\[
v=\langle z,z\rangle +\sum_{s,t\geq 2}F_{st}^{*}\left( z,\overline{z}%
,u\right) .
\]
Then the function $F_{32}^{*}(z,\overline{z},u)$ is given by 
\begin{eqnarray*}
F_{32}^{*}(z,\overline{z},u) \\
=4uA_{\alpha \beta \gamma \delta \overline{\eta }\overline{\xi }}z^{\alpha
}z^{\beta }z^{\gamma }a^{\delta }z^{\overline{\eta }}z^{\overline{\xi }%
}+3uA_{\alpha \beta \gamma \overline{\delta }\overline{\eta }\overline{\xi }%
}z^{\alpha }z^{\beta }z^{\gamma }z^{\overline{\delta }}z^{\overline{\eta }%
}a^{\overline{\xi }} \\
-\frac{6u\langle z,z\rangle ^{2}}{(n+1)(n+2)}\{4A_{\alpha \beta \gamma
\delta ..}^{\text{ }\alpha \beta }z^{\gamma }a^{\delta }+3A_{\alpha \beta
\gamma \overline{\delta }..}^{\text{ }\alpha \beta }z^{\gamma }a^{\overline{%
\delta }}\}.
\end{eqnarray*}
\end{lemma}

\proof
Suppose that $M$ is given up to weight $7$ by the following equation: 
\begin{eqnarray*}
v &=&\langle z,z\rangle +A_{\alpha \beta \gamma \delta \overline{\eta }%
\overline{\xi }}z^{\alpha }z^{\beta }z^{\gamma }z^{\delta }z^{\overline{\eta 
}}z^{\overline{\xi }}+A_{\alpha \beta \gamma \overline{\delta }\overline{%
\eta }\overline{\xi }}z^{\alpha }z^{\beta }z^{\gamma }z^{\overline{\delta }%
}z^{\overline{\eta }}z^{\overline{\xi }} \\
&&+A_{\alpha \beta \overline{\gamma }\overline{\delta }\overline{\eta }%
\overline{\xi }}z^{\alpha }z^{\beta }z^{\overline{\gamma }}z^{\overline{%
\delta }}z^{\overline{\eta }}z^{\overline{\xi }} \\
&&+C_{\alpha \beta \gamma \delta \eta \overline{\xi }\overline{\rho }%
}z^{\alpha }z^{\beta }z^{\gamma }z^{\delta }z^{\eta }z^{\overline{\xi }}z^{%
\overline{\rho }}+C_{\alpha \beta \gamma \delta \overline{\eta }\overline{%
\xi }\overline{\rho }}z^{\alpha }z^{\beta }z^{\gamma }z^{\delta }z^{%
\overline{\eta }}z^{\overline{\xi }}z^{\overline{\rho }} \\
&&+C_{\alpha \beta \gamma \overline{\delta }\overline{\eta }\overline{\xi }%
\overline{\rho }}z^{\alpha }z^{\beta }z^{\gamma }z^{\overline{\delta }}z^{%
\overline{\eta }}z^{\overline{\xi }}z^{\overline{\rho }}+C_{\alpha \beta 
\overline{\gamma }\overline{\delta }\overline{\eta }\overline{\xi }\overline{%
\rho }}z^{\alpha }z^{\beta }z^{\overline{\gamma }}z^{\overline{\delta }}z^{%
\overline{\eta }}z^{\overline{\xi }}z^{\overline{\rho }} \\
&&+\sum_{\min (s,t)\geq 2,s+t=8}O(z^{s}\overline{z}^{t})+\sum_{\min
(s,t)\geq 2,s+t=6}O(z^{s}\overline{z}^{t}u)+O(9),
\end{eqnarray*}
where 
\begin{gather*}
\overline{A_{\alpha \beta \gamma \delta \overline{\eta }\overline{\xi }}}%
=A_{\eta \xi \overline{\alpha }\overline{\beta }\overline{\gamma }\overline{%
\delta }},\text{ }\overline{A_{\alpha \beta \gamma \overline{\delta }%
\overline{\eta }\overline{\xi }}}=A_{\delta \eta \xi \overline{\alpha }%
\overline{\beta }\overline{\gamma }}, \\
\overline{C_{\alpha \beta \gamma \delta \eta \overline{\xi }\overline{\rho }}%
}=C_{\xi \rho \overline{\alpha }\overline{\beta }\overline{\gamma }\overline{%
\delta }\overline{\eta }},\text{ }\overline{C_{\alpha \beta \gamma \delta 
\overline{\eta }\overline{\xi }\overline{\rho }}}=C_{\eta \xi \rho \overline{%
\alpha }\overline{\beta }\overline{\gamma }\overline{\delta }}, \\
A_{\alpha \beta \gamma ...}^{\text{ }\alpha \beta \gamma }=0,
\end{gather*}
and all barred and unbarred indices are respectively symmetric.

By the mapping $\psi $ in the decomposition of $N_{\sigma }=\varphi \circ
E\circ \psi ,$ $M$ is transformed up to weight $7$ to a real hypersurface as
follows: 
\begin{eqnarray*}
v &=&\langle z,z\rangle +F_{04}(z,\overline{z},u)+F_{40}(z,\overline{z},u) \\
&&+F_{13}(z,\overline{z},u)+F_{31}(z,\overline{z},u)+F_{22}(z,\overline{z},u)
\\
&&+F_{14}(z,\overline{z},u)+F_{41}(z,\overline{z},u)+F_{23}(z,\overline{z}%
,u)+F_{32}(z,\overline{z},u) \\
&&+F_{15}(z,\overline{z},u)+F_{51}(z,\overline{z},u)+F_{33}(z,\overline{z},u)
\\
&&+O(z^{2}\overline{z}^{4})+O(z^{4}\overline{z}^{2})+\sum_{\min (s,t)\geq
2,s+t=7}O(z^{s}\overline{z}^{t})+O(9),
\end{eqnarray*}
where 
\begin{eqnarray*}
F_{40}(z,\overline{z},u) &=&2u^{2}A_{\alpha \beta \gamma \delta \overline{%
\eta }\overline{\xi }}z^{\alpha }z^{\beta }z^{\gamma }z^{\delta }a^{%
\overline{\eta }}a^{\overline{\xi }} \\
F_{31}(z,\overline{z},u) &=&8u^{2}A_{\alpha \beta \gamma \delta \overline{%
\eta }\overline{\xi }}z^{\alpha }z^{\beta }z^{\gamma }a^{\delta }z^{%
\overline{\eta }}a^{\overline{\xi }}+6u^{2}A_{\alpha \beta \gamma \overline{%
\delta }\overline{\eta }\overline{\xi }}z^{\alpha }z^{\beta }z^{\gamma }z^{%
\overline{\delta }}a^{\overline{\eta }}a^{\overline{\xi }} \\
F_{41}(z,\overline{z},u) &=&2uA_{\alpha \beta \gamma \delta \overline{\eta }%
\overline{\xi }}z^{\alpha }z^{\beta }z^{\gamma }z^{\delta }z^{\overline{\eta 
}}a^{\overline{\xi }} \\
F_{51}(z,\overline{z},u) &=&-4iu\langle z,z\rangle A_{\alpha \beta \gamma
\delta \overline{\eta }\overline{\xi }}z^{\alpha }z^{\beta }z^{\gamma
}z^{\delta }a^{\overline{\eta }}a^{\overline{\xi }}+2uC_{\alpha \beta \gamma
\delta \eta \overline{\xi }\overline{\rho }}z^{\alpha }z^{\beta }z^{\gamma
}z^{\delta }z^{\eta }z^{\overline{\xi }}a^{\overline{\rho }} \\
&&+12iu\langle z,a\rangle A_{\alpha \beta \gamma \delta \overline{\eta }%
\overline{\xi }}z^{\alpha }z^{\beta }z^{\gamma }z^{\delta }z^{\overline{\eta 
}}a^{\overline{\xi }} \\
F_{22}(z,\overline{z},u) &=&12u^{2}\{A_{\alpha \beta \gamma \delta \overline{%
\eta }\overline{\xi }}z^{\alpha }z^{\beta }a^{\gamma }a^{\delta }z^{%
\overline{\eta }}z^{\overline{\xi }}+A_{\alpha \beta \overline{\gamma }%
\overline{\delta }\overline{\eta }\overline{\xi }}z^{\alpha }z^{\beta }z^{%
\overline{\gamma }}z^{\overline{\delta }}a^{\overline{\eta }}a^{\overline{%
\xi }}\} \\
&&+9u^{2}A_{\alpha \beta \gamma \overline{\delta }\overline{\eta }\overline{%
\xi }}z^{\alpha }z^{\beta }a^{\gamma }z^{\overline{\delta }}z^{\overline{%
\eta }}a^{\overline{\xi }} \\
F_{32}(z,\overline{z},u) &=&4uA_{\alpha \beta \gamma \delta \overline{\eta }%
\overline{\xi }}z^{\alpha }z^{\beta }z^{\gamma }a^{\delta }z^{\overline{\eta 
}}z^{\overline{\xi }}+3uA_{\alpha \beta \gamma \overline{\delta }\overline{%
\eta }\overline{\xi }}z^{\alpha }z^{\beta }z^{\gamma }z^{\overline{\delta }%
}z^{\overline{\eta }}a^{\overline{\xi }} \\
F_{33}(z,\overline{z},u) &=&24iu\langle z,z\rangle \{A_{\alpha \beta \gamma
\delta \overline{\eta }\overline{\xi }}z^{\alpha }z^{\beta }a^{\gamma
}a^{\delta }z^{\overline{\eta }}z^{\overline{\xi }}-A_{\alpha \beta 
\overline{\gamma }\overline{\delta }\overline{\eta }\overline{\xi }%
}z^{\alpha }z^{\beta }z^{\overline{\gamma }}z^{\overline{\delta }}a^{%
\overline{\eta }}a^{\overline{\xi }}\} \\
&&+8iu\left\{ \langle z,a\rangle A_{\alpha \beta \overline{\gamma }\overline{%
\delta }\overline{\eta }\overline{\xi }}z^{\alpha }z^{\beta }z^{\overline{%
\gamma }}z^{\overline{\delta }}z^{\overline{\eta }}a^{\overline{\xi }%
}-\right. \\
&&\hspace{3cm}\left. \langle a,z\rangle A_{\alpha \beta \gamma \delta 
\overline{\eta }\overline{\xi }}z^{\alpha }z^{\beta }z^{\gamma }a^{\delta
}z^{\overline{\eta }}z^{\overline{\xi }}\right\} \\
&&+A_{\alpha \beta \gamma \overline{\delta }\overline{\eta }\overline{\xi }%
}z^{\alpha }z^{\beta }z^{\gamma }z^{\overline{\delta }}z^{\overline{\eta }%
}z^{\overline{\xi }} \\
&&+12iu\left\{ \langle z,a\rangle A_{\alpha \beta \gamma \overline{\delta }%
\overline{\eta }\overline{\xi }}z^{\alpha }z^{\beta }a^{\gamma }z^{\overline{%
\delta }}z^{\overline{\eta }}z^{\overline{\xi }}-\right. \\
&&\hspace{3cm}\left. \langle a,z\rangle A_{\alpha \beta \gamma \overline{%
\delta }\overline{\eta }\overline{\xi }}z^{\alpha }z^{\beta }z^{\gamma }z^{%
\overline{\delta }}z^{\overline{\eta }}a^{\overline{\xi }}\right\} \\
&&+4u\{C_{\alpha \beta \gamma \delta \overline{\eta }\overline{\xi }%
\overline{\rho }}z^{\alpha }z^{\beta }z^{\gamma }a^{\delta }z^{\overline{%
\eta }}z^{\overline{\xi }}z^{\overline{\rho }}+C_{\alpha \beta \gamma 
\overline{\delta }\overline{\eta }\overline{\xi }\overline{\rho }}z^{\alpha
}z^{\beta }z^{\gamma }z^{\overline{\delta }}z^{\overline{\eta }}z^{\overline{%
\xi }}a^{\overline{\rho }}\}.
\end{eqnarray*}

Since $M$ is in normal form, $E=id+O_{\times }(7)$ by Theorem \ref{Cor.3.2}.
So the function $p(u)$ satisfies 
\[
p(u)=\frac{1}{6}p^{\prime \prime \prime }(0)u^{3}+O(8). 
\]
Let $g(z,w)$ be a holomorphic function(cf. Lemma 3 in the paper \cite{Pa1})
defined by 
\begin{eqnarray*}
g(z,w)-g(0,w) &=&-2iF(p(w),,\overline{p}(w),w) \\
&&+2iF\left( z+p(w),\overline{p}(w),w+\frac{1}{2}\{g(z,w)-g(0,w)\}\right) \\
g(0,w) &=&iF(p(w),,\overline{p}(w),w).
\end{eqnarray*}
Thus the function $g(0,u)$ satisfies 
\[
g(0,u)=iF(p(u),\overline{p}(u),u)=O(12). 
\]
Hence the holomorphic function $g(z,w)$ is implicitly defined up to weight $%
11$ as follows: 
\[
g=2iF\left( z+p(w),\overline{p}(w),w+\frac{1}{2}g\right) +O(12). 
\]
Thus we obtain 
\begin{eqnarray*}
g &=&2iF\left( z+p(w),\overline{p}(w),w+\frac{1}{2}g\right) +O(12) \\
&=&\frac{i}{3}\langle z,p^{\prime \prime \prime }(0)\rangle w^{3}+2iF\left(
z,0,w+\frac{1}{2}g\right) +O(9) \\
&=&\frac{i}{3}\langle z,p^{\prime \prime \prime }(0)\rangle w^{3}+4i(w+\frac{%
1}{2}g)^{2}A_{\alpha \beta \gamma \delta \overline{\eta }\overline{\xi }%
}z^{\alpha }z^{\beta }z^{\gamma }z^{\delta }a^{\overline{\eta }}a^{\overline{%
\xi }}+O(9) \\
&=&\frac{i}{3}\langle z,p^{\prime \prime \prime }(0)\rangle
w^{3}+4iw^{2}A_{\alpha \beta \gamma \delta \overline{\eta }\overline{\xi }%
}z^{\alpha }z^{\beta }z^{\gamma }z^{\delta }a^{\overline{\eta }}a^{\overline{%
\xi }}+O(9).
\end{eqnarray*}

We carry out the following mapping: 
\begin{eqnarray*}
z &=&z^{*}+\frac{1}{6}p^{\prime \prime \prime }(0)w^{*3}+O(8) \\
w &=&w^{*}+\frac{i}{3}\langle z^{*},p^{\prime \prime \prime }(0)\rangle
w^{*3}+4iw^{*2}A_{\alpha \beta \gamma \delta \overline{\eta }\overline{\xi }%
}z^{*\alpha }z^{*\beta }z^{*\gamma }z^{*\delta }a^{\overline{\eta }}a^{%
\overline{\xi }}+O(9).
\end{eqnarray*}
Its inverse mapping is given as follows: 
\begin{align}
z^{*}=& z-\frac{1}{6}p^{\prime \prime \prime }(0)w^{3}+O(8)  \nonumber \\
w^{*}=& w-\frac{i}{3}\langle z,p^{\prime \prime \prime }(0)\rangle
w^{3}-4iw^{2}A_{\alpha \beta \gamma \delta \overline{\eta }\overline{\xi }%
}z^{\alpha }z^{\beta }z^{\gamma }z^{\delta }a^{\overline{\eta }}a^{\overline{%
\xi }}+O(9).  \tag*{(3.2)}  \label{car}
\end{align}
Then $M$ is transformed by the mapping \ref{car} to a real hypersurface up
to weight $6$ as follows: 
\begin{eqnarray*}
v &=&\langle z,z\rangle -iu^{2}\langle z,z\rangle \langle z,p^{\prime \prime
\prime }(0)\rangle +iu^{2}\langle z,z\rangle \langle p^{\prime \prime \prime
}(0),z\rangle \\
&&-4iu\langle z,z\rangle A_{\alpha \beta \gamma \delta \overline{\eta }%
\overline{\xi }}z^{\alpha }z^{\beta }z^{\gamma }z^{\delta }a^{\overline{\eta 
}}a^{\overline{\xi }}+4iu\langle z,z\rangle A_{\alpha \beta \overline{\gamma 
}\overline{\delta }\overline{\eta }\overline{\xi }}a^{\alpha }a^{\beta }z^{%
\overline{\gamma }}z^{\overline{\delta }}z^{\overline{\eta }}z^{\overline{%
\xi }} \\
&&+F_{13}(z,\overline{z},u)+F_{31}(z,\overline{z},u)+F_{22}(z,\overline{z},u)
\\
&&+F_{14}(z,\overline{z},u)+F_{41}(z,\overline{z},u)+F_{23}(z,\overline{z}%
,u)+F_{32}(z,\overline{z},u) \\
&&+F_{15}(z,\overline{z},u)+F_{51}(z,\overline{z},u)+F_{33}(z,\overline{z},u)
\\
&&+O(z^{2}\overline{z}^{4})+O(z^{4}\overline{z}^{2})+\sum_{\min (s,t)\geq
2,s+t=7}O(z^{s}\overline{z}^{t})+O(9) \\
&=&\langle z,z\rangle +F_{12}(z,\overline{z},u)+F_{21}(z,\overline{z}%
,u)+F_{13}(z,\overline{z},u)+F_{31}(z,\overline{z},u) \\
&&+F_{22}(z,\overline{z},u)+F_{14}(z,\overline{z},u)+F_{41}(z,\overline{z}%
,u)+F_{23}(z,\overline{z},u) \\
&&+F_{32}(z,\overline{z},u)+F_{15}^{*}(z,\overline{z},u)+F_{51}^{*}(z,%
\overline{z},u)+F_{33}(z,\overline{z},u) \\
&&+O(z^{2}\overline{z}^{4})+O(z^{4}\overline{z}^{2})+\sum_{\min (s,t)\geq
2,s+t=7}O(z^{s}\overline{z}^{t})+O(9),
\end{eqnarray*}
where 
\begin{eqnarray*}
F_{21}^{*}(z,\overline{z},u) &=&-iu^{2}\langle z,z\rangle \langle
z,p^{\prime \prime \prime }(0)\rangle \\
F_{51}^{*}(z,\overline{z},u) &=&F_{51}(z,\overline{z},u)-4iu\langle
z,z\rangle A_{\alpha \beta \gamma \delta \overline{\eta }\overline{\xi }%
}z^{\alpha }z^{\beta }z^{\gamma }z^{\delta }a^{\overline{\eta }}a^{\overline{%
\xi }}.
\end{eqnarray*}

We carry out the following mapping: 
\begin{eqnarray*}
z^{*\alpha } &=&z^{\alpha }-iz^{\alpha }\langle z,p^{\prime \prime \prime
}(0)\rangle w^{2}+8g^{\alpha \overline{\xi }}A_{\beta \gamma \delta \eta 
\overline{\xi }\overline{\rho }}z^{\beta }z^{\gamma }z^{\delta }a^{\eta }a^{%
\overline{\rho }}w^{2} \\
&&+6g^{\alpha \overline{\xi }}A_{\beta \gamma \delta \overline{\eta }%
\overline{\xi }\overline{\rho }}z^{\beta }z^{\gamma }z^{\delta }a^{\overline{%
\eta }}a^{\overline{\rho }}w^{2} \\
&&+2g^{\alpha \overline{\xi }}A_{\beta \gamma \delta \eta \overline{\xi }%
\overline{\rho }}z^{\beta }z^{\gamma }z^{\delta }z^{\eta }a^{\overline{\rho }%
}w-8iz^{\alpha }A_{\beta \gamma \delta \eta \overline{\xi }\overline{\rho }%
}z^{\beta }z^{\gamma }z^{\delta }z^{\eta }a^{\overline{\xi }}a^{\overline{%
\rho }}w \\
&&+12i\langle z,a\rangle g^{\alpha \overline{\xi }}A_{\beta \gamma \delta
\eta \overline{\xi }\overline{\rho }}z^{\beta }z^{\gamma }z^{\delta }z^{\eta
}a^{\overline{\rho }}w+2g^{\alpha \overline{\rho }}C_{\beta \gamma \delta
\eta \xi \overline{\rho }\overline{\sigma }}z^{\beta }z^{\gamma }z^{\delta
}z^{\eta }z^{\xi }a^{\overline{\sigma }}w \\
&&+O(8) \\
w^{*} &=&w.
\end{eqnarray*}
Then the real hypersurface is transformed to 
\begin{eqnarray*}
v &=&\langle z,z\rangle -2u\langle z,z\rangle ^{2}\langle z,p^{\prime \prime
\prime }(0)\rangle -2u\langle z,z\rangle ^{2}\langle p^{\prime \prime \prime
}(0),z\rangle \\
&&+F_{22}(z,\overline{z},u)+F_{23}(z,\overline{z},u)+F_{32}(z,\overline{z}%
,u)+F_{33}(z,\overline{z},u) \\
&&+O(z^{2}\overline{z}^{4})+O(z^{4}\overline{z}^{2})+\sum_{\min (s,t)\geq
2,s+t=7}O(z^{s}\overline{z}^{t})+O(9) \\
&=&\langle z,z\rangle +F_{22}(z,\overline{z},u)+F_{23}^{*}(z,\overline{z}%
,u)+F_{32}^{*}(z,\overline{z},u)+F_{33}(z,\overline{z},u) \\
&&+O(z^{2}\overline{z}^{4})+O(z^{4}\overline{z}^{2})+\sum_{\min (s,t)\geq
2,s+t=7}O(z^{s}\overline{z}^{t})+O(9),
\end{eqnarray*}
where 
\begin{align}
F_{32}^{*}(z,\overline{z},u)=& F_{32}(z,\overline{z},u)-2u\langle z,z\rangle
^{2}\langle z,p^{\prime \prime \prime }(0)\rangle  \nonumber \\
=& 4uA_{\alpha \beta \gamma \delta \overline{\eta }\overline{\xi }}z^{\alpha
}z^{\beta }z^{\gamma }a^{\delta }z^{\overline{\eta }}z^{\overline{\xi }%
}+3uA_{\alpha \beta \gamma \overline{\delta }\overline{\eta }\overline{\xi }%
}z^{\alpha }z^{\beta }z^{\gamma }z^{\overline{\delta }}z^{\overline{\eta }%
}a^{\overline{\xi }}  \nonumber \\
& -2u\langle z,z\rangle ^{2}\langle z,p^{\prime \prime \prime }(0)\rangle . 
\tag*{(3.3)}  \label{F32}
\end{align}
The condition $\Delta ^{2}F_{23}^{*}=0$ determines the coefficient $%
p^{\prime \prime \prime }(0)$ so that 
\begin{equation}
\langle z,p^{\prime \prime \prime }(0)\rangle =\frac{3}{(n+1)(n+2)}%
\{4A_{\alpha \beta \gamma \delta ..}^{\text{ }\alpha \beta }z^{\gamma
}a^{\delta }+3A_{\alpha \beta \gamma \overline{\delta }..}^{\text{ }\alpha
\beta }z^{\gamma }a^{\overline{\delta }}\}.  \tag*{(3.4)}  \label{F322}
\end{equation}

Since $E=id+O_{\times }(7),$ we consider a matrix $B$ such that 
\[
\langle z,z\rangle =\langle z+u^{3}Bz,z+u^{3}Bz\rangle +O(9), 
\]
which is equivalent to the following condition: 
\[
\langle Bz,z\rangle +\langle z,Bz\rangle =0. 
\]
Hence, on the real hypersurface, we have 
\begin{eqnarray*}
\langle z,z\rangle &=&\langle z+u^{3}Bz,z+u^{3}Bz\rangle +O(9) \\
&=&\langle z+w^{3}Bz,z+w^{3}Bz\rangle -6iu^{2}\langle z,z\rangle \langle
Bz,z\rangle +2i\langle z,z\rangle ^{3}\langle Bz,z\rangle \\
&&+O(9).
\end{eqnarray*}

We carry out the mapping: 
\begin{eqnarray*}
z^{*} &=&z+Bzw^{3} \\
w^{*} &=&w,
\end{eqnarray*}
Then the real hypersurface $M$ is transformed up to weight $8$ to 
\begin{eqnarray*}
v &=&\langle z,z\rangle +F_{22}^{*}(z,\overline{z},u)+F_{23}^{*}(z,\overline{%
z},u)+F_{32}^{*}(z,\overline{z},u)+F_{33}(z,\overline{z},u) \\
&&+O(z^{2}\overline{z}^{4})+O(z^{4}\overline{z}^{2})+\sum_{\min (s,t)\geq
2,s+t=7}O(z^{s}\overline{z}^{t})+O(9),
\end{eqnarray*}
where 
\[
F_{22}^{*}(z,\overline{z},u)=F_{22}(z,\overline{z},u)-6iu^{2}\langle
z,z\rangle \langle Bz,z\rangle . 
\]

The condition $\Delta ^{2}F_{22}^{*}=0$ determines the matrix $B$ so that 
\begin{eqnarray*}
\langle Bz,z\rangle &=&-\frac{2i}{n+2}\left\{ 4A_{\rho \beta \gamma \delta .%
\overline{\xi }}^{\text{ }\rho }z^{\beta }a^{\gamma }a^{\delta }z^{\overline{%
\xi }}+4A_{\rho \beta .\overline{\delta }\overline{\eta }\overline{\xi }}^{%
\text{ }\rho }z^{\beta }z^{\overline{\delta }}a^{\overline{\eta }}a^{%
\overline{\xi }}+\right. \\
&&\hspace{5cm}\left. 3A_{\rho \beta \gamma .\overline{\eta }\overline{\xi }%
}^{\text{ }\rho }z^{\beta }a^{\gamma }z^{\overline{\eta }}a^{\overline{\xi }%
}\right\} \\
&&+\frac{i\langle z,z\rangle }{(n+1)(n+2)}\left\{ 4A_{\rho \sigma \gamma
\delta ..}^{\text{ }\rho \sigma }a^{\gamma }a^{\delta }+4A_{\rho \sigma
\beta ..\overline{\eta }\overline{\xi }}^{\text{ }\rho \sigma }a^{\overline{%
\eta }}a^{\overline{\xi }}+\right. \\
&&\hspace{5cm}\left. 3A_{\rho \sigma \gamma ..\overline{\xi }}^{\text{ }\rho
\sigma }a^{\gamma }a^{\overline{\xi }}\right\} .
\end{eqnarray*}
Then we obtain $\Delta ^{3}F_{33}^{*}=0$ by putting 
\[
F_{33}^{*}(z,\overline{z},u)=F_{33}(z,\overline{z},u)-\kappa u\langle
z,z\rangle ^{3} 
\]
where 
\begin{eqnarray*}
\kappa &=&\frac{48i}{n(n+2)}\left\{ A_{\rho \sigma \gamma \delta ..}^{\text{ 
}\rho \sigma }a^{\gamma }a^{\delta }-A_{\rho \sigma ..\overline{\eta }%
\overline{\xi }}^{\text{ }\rho \sigma }a^{\overline{\eta }}a^{\overline{\xi }%
}\right\} \\
&&+\frac{24}{n(n+1)(n+2)}\left\{ C_{\alpha \beta \gamma \delta ...}^{\text{ }%
\alpha \beta \gamma }a^{\delta }+C_{\alpha \beta \gamma ...\overline{\rho }%
}^{\text{ }\alpha \beta \gamma }a^{\overline{\rho }}\right\} .
\end{eqnarray*}
Thus we have completed the normalizing process up to weight $8.$ Therefore,
the desired result is obtained by the equalities \ref{F32} and \ref{F322}.
This completes the proof.\endproof

\textbf{II}. Main theorem

\begin{theorem}
\label{Thm.5.18}Let $M$ be a nondegenerate analytic real hypersurface in a
complex manifold and $U$ be an open subset of $M$ consisting of umbilic
points. Then the open subset $U$ is locally biholomorphic to a real
hyperquadric.
\end{theorem}

\proof
For the case of $n=1,$ we have $F_{22}=F_{23}=F_{33}=0$. By the definition
on umbilic points for $n=1$ and invariance of normal form under translation
along $u$-curve, in any normal coordinate at any point of $U$ we have 
\[
F_{24}=F_{42}=0. 
\]
Suppose that there is a positive integer $k\geq 7$ and a normal coordinate
at a point $p$ in $U$ such that 
\begin{equation}
v=z\bar{z}+\sum_{\min (s,t)\geq 2,s+t=k}F_{st}(z,\bar{z},u)+\sum_{\min
(s,t)\geq 2,s+t\geq k+1}F_{st}(z,\bar{z},u),  \tag*{(3.5)}  \label{4.70}
\end{equation}
where 
\[
\sum_{\min (s,t)\geq 2,s+t=k}F_{st}(z,\bar{z},u)\neq 0. 
\]
Let's take a point $p\in U$ and a normal coordinate so that $k$ is the
smallest integer satisfying the condition in \ref{4.70}. Since the normal
form is invariant under translation along the $u$-curve, we can assume that 
\[
\sum_{\min (s,t)\geq 2,s+t=k}F_{st}(z,\bar{z},0)\neq 0. 
\]
By Lemma \ref{Lem.5.19} this is a contradiction to the choice of the integer 
$k$. Therefore $M$ is defined by 
\[
v=z\bar{z} 
\]
in any normal coordinate at any point in $U$. This completes the proof for $%
n=1$.

For the case of $n\geq 2$, by Lemma \ref{Lem.5.19} it suffices to show that 
\[
F_{22}=F_{23}=F_{24}=F_{33}=0. 
\]
Since the normal form is invariant under translation along $u$-curve, we
have 
\[
F_{22}=0 
\]
in any normal coordinate at any point in $U$.

Suppose that there is a point $p\in U$ and a normal coordinate at $p$ such
that $F_{23}\neq 0$. Further, without loss of generality we can assume that $%
F_{23}(z,\bar{z},0)\neq 0$ so that 
\[
F_{23}(z,\bar{z},0)=C_{\alpha \beta \bar{\gamma}\bar{\delta}\bar{\eta}%
}z^{\alpha }z^{\beta }z^{\bar{\gamma}}z^{\bar{\delta}}z^{\bar{\eta}}\neq 0. 
\]
Then from the proof of Theorem \ref{Th5} we obtain the following identity
for all $a\in \Bbb{C}^{n}$%
\begin{eqnarray*}
F_{22}^{*}(z,\overline{z},u) &=&3uC_{\alpha \beta \gamma \bar{\delta}%
\overline{\eta }}z^{\alpha }z^{\beta }a^{\gamma }z^{\bar{\delta}}z^{%
\overline{\eta }}+3uC_{\alpha \beta \overline{\gamma }\bar{\delta}\overline{%
\eta }}z^{\alpha }z^{\beta }z^{\overline{\gamma }}z^{\bar{\delta}}a^{%
\overline{\eta }} \\
&&-\frac{12}{n+2}\cdot u\langle z,z\rangle \left\{ C_{\delta \alpha \beta .%
\overline{\gamma }}^{\text{ }\delta }z^{\alpha }a^{\beta }z^{\overline{%
\gamma }}+C_{\delta \alpha .\overline{\beta }\overline{\gamma }}^{\text{ }%
\delta }z^{\alpha }z^{\overline{\beta }}a^{\overline{\gamma }}\right\} \\
&=&0.
\end{eqnarray*}
Hence we obtain 
\[
(n+2)C_{\alpha \beta \gamma \bar{\delta}\overline{\eta }}=h_{\alpha 
\overline{\delta }}C_{\rho \beta \gamma .\overline{\eta }}^{\text{ }\rho
}+h_{\beta \overline{\delta }}C_{\rho \alpha \gamma .\overline{\eta }}^{%
\text{ }\rho }+h_{\alpha \overline{\eta }}C_{\rho \beta \gamma .\overline{%
\delta }}^{\text{ }\rho }+h_{\beta \overline{\eta }}C_{\rho \alpha \gamma .%
\overline{\delta }}^{\text{ }\rho }. 
\]
Contracting the pair $(\gamma ,\overline{\delta })$ yields 
\[
C_{\rho \beta \gamma .\overline{\eta }}^{\text{ }\rho }=0, 
\]
which implies 
\[
C_{\alpha \beta \gamma \bar{\delta}\overline{\eta }}=0. 
\]
This is a contradiction to the hypothesis that $F_{23}(z,\bar{z}%
,0)=C_{\alpha \beta \bar{\gamma}\bar{\delta}\bar{\eta}}z^{\alpha }z^{\beta
}z^{\bar{\gamma}}z^{\bar{\delta}}z^{\bar{\eta}}\neq 0.$ Thus 
\[
F_{23}=F_{32}=0 
\]
in any normal coordinate at any point in $U$.

Suppose that there is a point $p\in U$ and a normal coordinate at the point
such that 
\[
F_{24}+F_{33}+F_{42}\neq 0. 
\]
Without loss of generality, we may assume 
\[
F_{24}(z,\bar{z},0)\neq 0\quad \text{or\quad }F_{33}(z,\bar{z},0)\neq 0. 
\]
Let's put 
\begin{eqnarray*}
F_{24}(z,\bar{z},0) &=&A_{\alpha \beta \bar{\gamma}\bar{\delta}\bar{\eta}%
\bar{\sigma}}z^{\alpha }z^{\beta }z^{\bar{\gamma}}z^{\bar{\delta}}z^{\bar{%
\eta}}z^{\bar{\sigma}} \\
F_{33}(z,\bar{z},0) &=&A_{\alpha \beta \gamma \bar{\delta}\bar{\eta}\bar{%
\sigma}}z^{\alpha }z^{\beta }z^{\gamma }z^{\bar{\delta}}z^{\bar{\eta}}z^{%
\bar{\sigma}}.
\end{eqnarray*}
Since $F_{32}^{*}(z,\bar{z},u)=0$ identically, Lemma \ref{new} yields the
following identity for all $a\in \Bbb{C}^{n}$: 
\begin{eqnarray*}
&&4A_{\alpha \beta \gamma \delta \overline{\eta }\overline{\xi }}z^{\alpha
}z^{\beta }z^{\gamma }a^{\delta }z^{\overline{\eta }}z^{\overline{\xi }%
}+3A_{\alpha \beta \gamma \overline{\delta }\overline{\eta }\overline{\xi }%
}z^{\alpha }z^{\beta }z^{\gamma }z^{\overline{\delta }}z^{\overline{\eta }%
}a^{\overline{\xi }} \\
&=&\frac{6\langle z,z\rangle ^{2}}{(n+1)(n+2)}\left\{ 4A_{\alpha \beta
\gamma \delta ..}^{\text{ }\alpha \beta }z^{\gamma }a^{\delta }+3A_{\alpha
\beta \gamma \overline{\delta }..}^{\text{ }\alpha \beta }z^{\gamma }a^{%
\overline{\delta }}\right\} .
\end{eqnarray*}
Hence we obtain 
\begin{eqnarray*}
A_{\alpha \beta \gamma \delta \overline{\eta }\overline{\xi }}z^{\alpha
}z^{\beta }z^{\gamma }a^{\delta }z^{\overline{\eta }}z^{\overline{\xi }} &=&%
\frac{6}{(n+1)(n+2)}\langle z,z\rangle ^{2}A_{\alpha \beta \gamma \delta
..}^{\text{ }\alpha \beta }z^{\gamma }a^{\delta }, \\
A_{\alpha \beta \gamma \overline{\delta }\overline{\eta }\overline{\xi }%
}z^{\alpha }z^{\beta }z^{\gamma }z^{\overline{\delta }}z^{\overline{\eta }%
}a^{\overline{\xi }} &=&\frac{6}{(n+1)(n+2)}\langle z,z\rangle ^{2}A_{\alpha
\beta \gamma ..\overline{\delta }}^{\text{ }\alpha \beta }z^{\gamma }a^{%
\overline{\delta }}.
\end{eqnarray*}
Putting $a^{\alpha }=z^{\alpha }$ and computing the derivative $\Delta ^{2}$
yields 
\[
A_{\alpha \beta \gamma \delta \overline{\eta }\overline{\xi }}=A_{\alpha
\beta \gamma \overline{\delta }\overline{\eta }\overline{\xi }}=0. 
\]
This is a contradiction to the hypothesis 
\[
F_{24}(z,\bar{z},0)\neq 0\quad \text{or\quad }F_{33}(z,\bar{z},0)\neq 0. 
\]
So $F_{24}=F_{42}=F_{33}=0$ in any normal coordinate at any point in $U$.

Then, by Lemma \ref{Lem.5.19}, $M$ is defined by 
\[
v=\langle z,z\rangle 
\]
in any normal coordinate at any point in $U$. This completes the proof.\endproof

There is a global version of Theorem \ref{CK}. We shall give a complete
proof of this global version to the paper \cite{Pa3}, where we present a new
proof of the following theorem:

\begin{lemma}[Pinchuk]
\label{CCM}Let $M$ be a nondegenerate analytic real hypersurface. Suppose
that there are two points $p,q$ on $M$ and a curve $\gamma $ on $M$
connecting the two points $p,q.$ Then $p$ is a spherical point if and only
if $q$ is a spherical point.
\end{lemma}

From Lemma \ref{CCM}, we obtain the following global result:

\begin{lemma}
\label{Tm8}Let $M$ be a connected nondegenerate analytic real hypersurface.
Suppose that the subset of umbilic points of $M$ has an interior point. Then 
$M$ is locally biholomorphic to a real hyperquadric at every point of $M.$
\end{lemma}

Then we obtain the following theorem from Theorem \ref{Thm.5.18} and Lemma 
\ref{Tm8}

\begin{theorem}
\label{Co8}Let $M$ be a connected nondegenerate analytic real hypersurface
in a complex manifold. Suppose that there is a point $p$ on $M$ at which $M$
is not locally biholomorphic to a real hyperquadric. Then $M$ is not locally
biholomorphic to a real hyperquadric at every point of $M.$
\end{theorem}

\textbf{III}. A family of nondegenerate analytic real hypersurfaces $M$
parametrized by $\mu .$ From Theorem \ref{main}, we obtain

\begin{lemma}
\label{Cor.3.4}Let $M_{\mu }$ be a nondegenerate analytic real hypersurface
defined by 
\[
v=F(z,\bar{z},u;\mu ),\quad \left. F\right| _{0}=\left. dF\right| _{0}=0,
\]
where $\mu $ is a real parameter and $F(z,\bar{z},u;\mu )$ is real analytic
of $\mu $ near $\mu =0.$ Let $\phi =(f,g)$ be a normalization of $M_{\mu }$
with initial value independent of $\mu .$ Suppose that $\phi (M_{\mu })$ is
defined by $v=\langle z,z\rangle +F^{*}\left( z,\bar{z},u\right) .$ Then the
functions $f,g,F^{*}$ are analytic of the parameter $\mu $ such that $%
f,g,F^{*}$ mod $\mu ^{l}$ is completely determined by $F$ mod $\mu ^{l}$ for
each $l\geq 0.$
\end{lemma}

By Theorem \ref{main}and Lemma \ref{Cor.3.4}, a normalization $\phi $ of $%
M_{\mu }$ and the operation of mod $\mu ^{l}$ are commutative, i.e., the
following two results give the same information up to $\rm{mod}$ $\mu ^{l}$
for each $l:$%
\[
\left\{ 
\begin{array}{l}
\phi \left( M_{\mu }\right) \quad \rm{mod}\mu ^{l} \\ 
\phi \left( M_{\mu }\quad \rm{mod}\mu ^{l}\right) 
\end{array}
\right. .
\]
Let $N_{\sigma }=\varphi \circ E\circ \psi $ be a normalization of $M$ and $%
M^{\prime }\equiv \psi \left( M\right) .$ Then $M^{\prime }$ depends
analytically on the parameter $a$ in the initial value $\sigma =(C,a,\rho
,r).$ By Lemma \ref{Cor.3.4} we can compute the normalization $E$ up to a
given order of $a$ inclusive.

\begin{lemma}
\label{Lem6}Let $k$ be a positive integer$\geq 7.$ If $M$ is an analytic
real hypersurface in normal form such that 
\[
v=\langle z,z\rangle +\sum_{\min (s,t)\geq 2,s+t=k}F_{st}(z,\bar{z}%
,u)+\sum_{\min (s,t)\geq 2,s+t\geq k+1}F_{st}(z,\bar{z},u),
\]
where 
\[
\sum_{\min (s,t)\geq 2,s+t=k}F_{st}(z,\bar{z},u)\neq 0.
\]
\ Then there is a vector $a\in \Bbb{C}^{n}$ and a normalization $N_{\sigma },
$ $\sigma =(id_{n\times n},a,1,0),$ such that 
\[
N_{\sigma }(M):v=\langle z,z\rangle +\sum_{\min (s,t)\geq
2}^{s+t=k-1}F_{st}^{*}(z,\bar{z},u)+\sum_{\min (s,t)\geq 2,s+t\geq
k}F_{st}^{*}(z,\bar{z},u),
\]
where 
\[
\sum_{\min (s,t)\geq 2,s+t=k-1}F_{st}^{*}(z,\bar{z},u)\neq 0.
\]
\end{lemma}

\proof
The real hypersurface $M$ is defined as follows: 
\[
v=\langle z,z\rangle +\sum_{\min (s,t)\geq 2,s+t=k}^{s+t=k+1}F_{st}(z,\bar{z}%
,u)+\sum_{s+t\geq k+2}O(z^{s}\bar{z}^{t}). 
\]
Let $M^{\prime }$ be a real hypersurface obtained after the mapping $\psi $
in the decomposition 
\[
N_{\sigma }=\varphi \circ E\circ \psi 
\]
where $\phi _{\sigma }=\varphi \circ \psi $ is the local automorphism of a
real hyperquadric with initial value $\sigma \in H$.

Then the real hypersurface $M^{\prime }$ depends analytically on the
parameter $a$ in $\sigma =(C,a,\rho ,r).$ By Lemma \ref{Cor.3.4} we can
compute the normalization $E$ up to $O(\left| a\right| ^{2}).$

By the mapping $\psi $ in the decomposition of $N_{\sigma }=\varphi \circ
E\circ \psi ,$ we obtain 
\begin{eqnarray*}
v &=&\langle z,z\rangle \\
&&+\sum_{\min (s,t)\geq 2,s+t=k}^{s+t=k+1}\sum_{\alpha }u\left\{ a^{\alpha
}\left( \frac{\partial F_{st}}{\partial z^{\alpha }}\right) (z,\bar{z},u)+%
\overline{a}^{\alpha }\left( \frac{\partial F_{st}}{\partial \overline{z}%
^{\alpha }}\right) (z,\bar{z},u)\right\} \\
&&+\sum_{\min (s,t)\geq 2,s+t=k}\sum_{\alpha }i\langle z,z\rangle \left\{
a^{\alpha }\left( \frac{\partial F_{st}}{\partial z^{\alpha }}\right) (z,%
\bar{z},u)-\overline{a}^{\alpha }\left( \frac{\partial F_{st}}{\partial 
\overline{z}^{\alpha }}\right) (z,\bar{z},u)\right\} \\
&&+\sum_{\min (s,t)\geq 2,s+t=k}\sum_{\alpha }i(\langle z,a\rangle -\langle
a,z\rangle )\left( \frac{\partial F_{st}}{\partial u}\right) (z,\bar{z},u) \\
&&+\sum_{\min (s,t)\geq 2,s+t=k}^{s+t=k}\left( 1-2i(1-s)\langle z,a\rangle
+2i(1-t)\langle a,z\rangle \right) F_{st}(z,\bar{z},u) \\
&&+\sum_{\min (s,t)\geq 2,s+t=k+2}O(z^{s}\bar{z}^{t})+O(\left| a\right|
^{2}).
\end{eqnarray*}
For $k\geq 8,$ by the normalization $E$ up to $O(\left| a\right| ^{2}),$ we
obtain 
\begin{eqnarray*}
v &=&\langle z,z\rangle \\
&&+\sum_{\min (s+1,t)\geq 2,s+t=k-1}\sum_{\alpha }ua^{\alpha }\left( \frac{%
\partial F_{s+1,t}}{\partial z^{\alpha }}\right) (z,\bar{z},u) \\
&&+\sum_{\min (s,t+1)\geq 2,s+t=k-1}\sum_{\alpha }u\overline{a}^{\alpha
}\left( \frac{\partial F_{s,t+1}}{\partial \overline{z}^{\alpha }}\right) (z,%
\bar{z},u) \\
&&+\sum_{\min (s,t)\geq 2,s+t=k}F_{st}(z,\bar{z},u)+\sum_{\min (s,t)\geq
2,s+t=k+1}O(z^{s}\bar{z}^{t})+O(\left| a\right| ^{2}).
\end{eqnarray*}
For $k=7,$ the condition $\Delta ^{3}F_{33}=0$ can be achieved without any
effects on the following term up to $O(\left| a\right| ^{2})$: 
\[
F_{42}^{*}(z,\overline{z},u)=\sum_{\alpha }\left\{ a^{\alpha }\left( \frac{%
\partial F_{52}}{\partial z^{\alpha }}\right) (z,\bar{z},u)+\overline{a}%
^{\alpha }\left( \frac{\partial F_{43}}{\partial \overline{z}^{\alpha }}%
\right) (z,\bar{z},u)\right\} +O(\left| a\right| ^{2}). 
\]
By the normalization $E$ up to $O(\left| a\right| ^{2}),$ we obtain 
\begin{eqnarray*}
v &=&\langle z,z\rangle \\
&&+\sum_{\min (s+1,t)\geq 2,s+t=6}\sum_{\alpha }ua^{\alpha }\left( \frac{%
\partial F_{s+1,t}}{\partial z^{\alpha }}\right) (z,\bar{z},u) \\
&&+\sum_{\min (s,t+1)\geq 2,s+t=6}\sum_{\alpha }u\overline{a}^{\alpha
}\left( \frac{\partial F_{s,t+1}}{\partial \overline{z}^{\alpha }}\right) (z,%
\bar{z},u) \\
&&+\kappa (u)\langle z,z\rangle ^{3}+\sum_{\min (s,t)\geq 2,s+t=7}F_{st}(z,%
\bar{z},u) \\
&&+\sum_{\min (s,t)\geq 2,s+t=8}O(z^{s}\bar{z}^{t})+O(\left| a\right| ^{2})
\end{eqnarray*}
where $\kappa (u)$ is determined so that $\Delta ^{3}F_{33}^{*}=0$ up to $%
O(\left| a\right| ^{2}).$

Since $N_{\sigma }(M)$ depends analytically on the parameter $a$ in $\sigma
=(C,a,\rho ,r)$(cf. \cite{Pa1})$,$ we obtain for sufficiently small $a\in 
\Bbb{C}^{n}$%
\[
\sum_{s+t=k-1}F_{st}^{*}(z,\overline{z},u)\neq 0. 
\]
This completes the proof.\endproof

\begin{theorem}
\label{last}Let $M$ be a nonspherical nondegenerate analytic real
hypersurface. Then, for $\dim M=3$, there exists $\sigma \in H$ such that 
\[
N_{\sigma }(M):v=\langle z,z\rangle +\sum_{\min (s,t)\geq 2,s+t\geq
6}F_{st}(z,\bar{z},u),
\]
where 
\[
F_{24}(z,\bar{z},u)\neq 0.
\]
For $\dim M\geq 5$, there exists $\sigma \in H$ such that 
\[
N_{\sigma }(M):v=\langle z,z\rangle +\sum_{\min (s,t)\geq 2}F_{st}(z,\bar{z}%
,u),
\]
where 
\[
F_{22}(z,\bar{z},u)\neq 0.
\]
\end{theorem}

We present a proof of Theorem \ref{last} in the paper \cite{Pa3}. Notice
that the case of $\dim M=3$ is an immediate consequence of Lemma \ref{Lem6}.

\end{document}